\documentclass[11pt]{article}
\usepackage[T1]{fontenc}
\usepackage[french]{babel}
\usepackage[latin1]{inputenc}
\usepackage{enumerate}
\usepackage{amsmath}
\usepackage{amssymb}
\usepackage{amsfonts}
\usepackage{mathrsfs}

\newcommand{\cp}{\mathbb{C}}
\newcommand{\pp}{\mathbb{P}}
\newcommand{\gp}{\mathbb{G}}
\newcommand{\om}{\omega}
\newcommand{\ep}{\epsilon}

\newcommand{\wt}{\widetilde}
\newtheorem{lem}{Lemme}

\newtheorem{prop}{Proposition}
\newtheorem{deff}{Définition}
\newtheorem{thr}{Th\'eor\`eme}
\newtheorem{rem}{Remarque}

\begin{document}
\date{27 Mai 2004}
\author{Martin Weimann\\{\footnotesize Laboratoire analyse et géométrie}\\{\footnotesize Université Bordeaux1}\\{\footnotesize 351, cours de la Libération}\\{\footnotesize 33045 TALENCE}}
\title{La trace {\it via} le calcul résiduel~: une nouvelle version du théorème d'Abel inverse, formes abéliennes}
\maketitle

{\footnotesize 
Le théorème d'Abel inverse affirme qu'il suffit qu'un germe de trace de forme m\'eromorphe 
soit rationnel pour qu'il se prolonge en la trace d'une forme rationnelle sur une sous-vari\'et\'e 
alg\'ebrique. Ce r\'esultat a été prouv\'e dans le cas trace nulle par 
P.A. Griffiths en 1976 \cite{G:gnus} puis dans le cas trace rationnelle par G. Henkin et M. Passare 
en 1999 \cite{hp:gnus}. On montre ici qu'il suffit que la trace 
soit rationnelle en les param\`etres qui ne correspondent pas aux pentes pour retrouver la conclusion du théorème d'Abel-inverse. 
La démonstration s'appuie essentiellement sur le calcul résiduel (courant résiduel, résidus ponctuels 
et théorème de dualité) qui permettent d'un côté de caractériser les formes traces d'une manière 
particulièrement algébrique et d'un autre côté de reconstituer un ensemble analytique et une forme 
méromorphe à partir des traces d'un nombre fini de fonctions. Cette démonstration permet d'établir le lien avec le théorème de Wood. De plus, on retrouve les bornes de 
Castelnuovo pour la dimension de l'espace des $q$-formes abéliennes sur une hypersurface de $\pp^{n+1}$. 
Les résultats présentés ici s'inspirent 
principalement des travaux d'Alain Yger \cite{Y:gnus}. Des idées similaires ont été développées récemment 
et indépendament par B. Fabre dans le cas plus général des courants localement résiduels \cite{F_2:gnus}.
}

\section{Introduction}

La donnée d'un sous-ensemble analytique ferm\'e $V$ de codimension pure $r$ d'un domaine $D\subset\pp^{n+1}(\cp)$ 
et d'une forme $(q,0)$ méromorphe $\Phi$ sur $V$ ($0\leq q\leq \dim V= 
n+1-r$) permettent de définir sur $D$ le courant $[V]\land \Phi$ à 
support dans $V$~; pour toute forme test $\phi$ (de bidegré $(n+1-r-q,n+1-r)$), on pose formellement~: 
$$
\langle [V]\land \Phi, \varphi\rangle=\int_{V} 
\Phi\land \varphi\,. 
$$
Plus pr\'ecis\'ement si $\delta$ d\'esigne le d\'enominateur universel 
pour les formes m\'eromorphes sur $V$ au voisinage d'un point $z\in V$, l'action du  
courant $[V]\wedge \Phi$ sur une forme de bidegr\'e $(n+1-r-q,n+1-r)$ de support dans un voisinage suffisamment petit 
de $z_0$ est d\'efinie comme la valeur en $\lambda=0$ du prolongement 
m\'eromorphe de la fonction 
$$
\lambda \to \int_{V\setminus {\rm Sing}(V)} |\delta|^{2\lambda} \Phi \wedge \varphi  
$$
(on montre que ce prolongement m\'eromorphe n'a pas de p\^ole en $\lambda=0$).  
On dit que $\Phi$ est régulière sur $V$ (au sens de Barlet, \cite{bar:gnus,hp:gnus}) 
si ce courant est $\bar{\partial}$-fermé dans $D$. 
Dans le cas où $\Phi$ est une forme de degré maximal sur $V$, ce courant est 
$d$-fermé en dehors du lieu polaire de $\Phi$. 
Si l'on suppose $D$ $r$-concave (c'est à dire réunion d'espaces linéaires de dimension $r$ appelés $r$-plans) 
avec $r=n+1-p={\rm codim}\, V$, on peut définir l'ouvert $D^*\subset\gp(r,n+1)$, $r$-dual de $D$, dont les éléments $t$ paramètrent les $r$-plans $L_t$ 
inclus dans $D$. On peut alors définir la variété d'incidence de $D$ comme étant 
la sous-vari\'et\'e analytique de $D\times D^*$ d\'efinie comme 
$$
{\rm INCID}_D:=\{(z,t)\in D\times D^* \;;\;z\in L_t\}.
$$
On note $p_1$ et $p_2$ les projections naturelles respectives de ${\rm INCID}_D$ sur $D$ et $D^*$. 
Le fait que $p_1$ soit une submersion et que $p_2$ soit propre sur le support du courant  $\Phi\land [V]$ 
(la propreté de $p_2$ est liée à 
la compacité de l'espace projectif) 
permettent de définir un nouveau courant appelé transformée d'Abel de $\Phi\land [V]$, 
noté $A(\Phi\land [V])$, ce en posant 
$$
A(\Phi\land [V]):=(p_2)_*(p_1^*(\Phi\land [V]))
$$
Ce nouveau courant est une forme méromorphe définie sur $D^*$ de même type que $\Phi$ 
que l'on appelle la trace  de $\Phi$ sur $V$ que l'on notera 
${\rm Tr}_V\, (\Phi)$ \footnote{en général la trace d'un courant désigne l'image directe de celui-ci par une application propre sur le support du courant. Ici c'est quelque peu différent car la trace correspond à l'image directe du courant $p_1^*(\Phi\land [V])$ par la projection sur $D^*$.}. Pour $t$ générique, le calcul de cette trace consiste 
à sommer les formes $\xi_j^* (\Phi)$ o\`u les points $\xi_j (t)$ sont les points d'intersections de 
$V$ avec $L_t$ (ces derniers étant en nombre fini, ce nombre restant localement constant). 
L'approche plus formelle {\it via} les courants a l'avantage de ne plus se soucier des singularités de $V$. 
On peut se réferer à \cite{hp:gnus} pour un rappel de ces diverses notions ainsi que pour une pr\'esentation 
de cette approche {\it via} la th\'eorie des courants. 
\vskip 2mm
\noindent
On s'intéresse ici au cas particulier où $V$ est une réunion finie de germes d'hypersurfaces 
(éventuellement singulières et non réduites) coupées (mais non incluses) 
par une droite de $D$ et $\Phi$ est un germe de forme méromorphe de degré maximal sur $V$ ({\it i.e.} 
de type $(n,0)$). Dans ce cas il s'avère que la donnée du couple $(V,\Phi)$ 
est équivalente à la donnée d'un couple $(F,H)$ de polynômes à une variable 
dont les coefficients sont des germes de fonctions méromorphes en les paramètres de 
$D^*$ vérifiant quelques conditions algébriques particulièrement simples. 
Bien que la démonstration de cette équivalence soit basée sur le théorème de préparation de 
Weierstrass, la construction des deux polynômes $F$ et $H$ se fait ais\'ement 
{\it via} le calcul résiduel. En effet, les coefficients de $F$ et $H$ sont uniquement déterminés 
par un système lin\'eaire non dégénéré dont les coefficients sont obtenus comme les traces 
d'un nombre fini de fonctions (ce sont donc des sommes compl\`etes de résidus). 
Ce phénomène permet de caractériser uniquement le couple $(V,\Phi)$ en terme des traces 
${\rm Tr}_V\, (y^k)$ et ${\rm Tr_V}\, (y^k\Phi)$, 
pour $k=0,..,2d-1$ (o\`u $d:=Tr_V\, 1$), \emph{calcul\'ees suivant 
une famille de droites d'une même direction donn\'ee} \footnote{en fait en restreignant à une seule projection, on retombe sur le concept de trace usuel originellement développé par Barlet \cite{bar:gnus}, et on retrouve ses résultats de caractérisation des formes régulières}. Ce nouveau résultat permet de 
voir la trace de $\Phi$ sur $V$ comme un calcul de résidus d'une fonction rationnelle d'une variable, 
c'est-à-dire comme un calcul de reste dans une simple division euclidienne (cette approche se trouvait 
d\'ej\`a en fait esquiss\'ee dans l'article de P.A. Griffiths \cite{G:gnus}). Ce point de vue permet de 
redémontrer le théorème d'Abel inverse sous une version plus forte~; en effet on montre alors que la 
rationnalité de la trace en $n$ des $2n$ variables de $D^*$ suffit pour conclure que $V$ est alors une 
courbe algébrique et $\Phi$ une forme rationnelle sur $V$. On comprendra ainsi le lien avec le théorème de Wood qui affirme qu'il est équivalent que $V$ soit algébrique de degré $d$ et que la trace de $y$ soit affine en $n$ des $2n$ variables (c'était en fait une des motivations de ce travail). Pour finir, on retrouve les bornes de Castelnuovo 
pour majorer la dimension de l'espace des $q$-formes abéliennes ({\it i.e.} de trace nulle), bornes 
atteintes dans le cas $q=n$ (voir \cite{he:gnus} par exemple pour une approche {\it via} les tissus).
\vskip 2mm
\noindent
L'aspect global du calcul résiduel permet de ne plus travailler localement avec 
le théorème des fonctions implicites et le lemme de Darboux (ingrédient essentiel 
dans la preuve d'Henkin et Passare) et donc de ne plus se soucier que les germes soient lisses ou réduits~;  
on esquive également ainsi les théorèmes d'Hartogs et de Remmert-Stein. 
De plus la preuve du théorème de préparation de Weierstrass peut s'obtenir à partir 
du calcul de résidus \cite{GH:gnus}~; il en est de m\^eme pour la construction du couple $(F,H)$  
et la trace est ainsi une forme différentielle dont les coefficients sont des sommes compl\`etes de 
résidus de polynômes à une variable, c'est-à-dire des restes de divisions euclidiennes 
dans l'anneau des polynômes. En ce sens, on peut parler de l'aspect très 
algébrique du théorème d'Abel et de son inversion, et plus 
généralement du concept de trace. Ce type de construction suggère de nouvelles démonstrations pour des 
généralisations du théorème d'Abel inverse, par exemple en remplaçant les espaces linéaires par une famille continuement paramétrée
d'ensembles algébriques de même dimension et de degré donné \cite{F_1:gnus}.

\section{Les coefficients de la trace exprim\'es en termes de sommes de r\'esidus}

Soit $V$ un sous-ensemble analytique de codimension $r$ d'un domaine $r$-concave $D\subset \pp^{n+1}(\cp)$. 
Soit $\Phi$ une $q$-forme méromorphe sur $V$ (c'est-à-dire localement la restriction à $V$ d'une forme 
méromorphe dans l'espace ambiant et dont le lieu polaire ${\rm Pol}\, (\Phi)$ v\'erifie 
$\dim ({\rm Pol}\, (\Phi) \cap V)<n+1-r=\dim V$). En utilisant l'interprétation de la trace 
en termes de courants, on montre ici par dualité que les coefficients de la trace 
s'expriment en termes de l'action de courants 
résiduels et sont en fait des sommes compl\`etes 
de résidus fonctions des paramètres $t$ de l'espace $r$-dual $D^*$ de $D$. 
Dans un premier temps on utilise les variables $Z$ pour $D$ et $t$ pour $D^*$. Pour toute forme-test $\varphi$  
(de type correctement choisi en fonction du type de $\Phi$) de $D^*$, on a par définition 
\begin{eqnarray*}
\langle {\rm Tr}_V\, (\Phi)\,,\,\varphi\rangle &=&
\int_{D^*} {\rm Tr}_V (\Phi) (t)\land \varphi(t) \\
&=&\int_{{\rm INCID}_D} ([p_1^{-1} (V)]\wedge 
p_1^* [\Phi]) \land {p_2}^* (\varphi)(Z,t)\,. 
\end{eqnarray*} 
En notant alors $[I_D]$ le courant d'intégration associé au sous-ensemble analytique 
${\rm INCID}_D$ de $D\times D^*$, on obtient l'expression
$$
\langle Tr_V\Phi,\varphi\rangle =\int_{D\times D^*} ([V](Z)\land[I_D](Z,t))\land \Phi(Z) 
\land \varphi(t)=\langle T\,,\, \varphi\rangle   
$$
où le $(q+n+1,n+1)$-courant $T:=([V(Z)]\land[I_D(Z,t)])\wedge \Phi(Z)$ a pour image 
directe par $p_2$ un courant forme méromorphe de même type que $\Phi$ correspondant pr\'ecis\'ement \`a la
trace de $\Phi$. On voudrait expliciter dans des cartes affines 
de la grassmannienne les coefficients de ce courant 
dans le cas d'une hypersurface $V\subset D$ (on retrouve alors la transformée d'Abel-Radon usuelle).
Quitte \`a faire agir un automorphisme de $\pp^{n+1} (\cp)$, on peut toujours supposer, si l'on 
note $(x_1,x_2,...,x_n,y,z)$ les coordonn\'ees homog\`enes dans $\pp^{n+1}(\cp)$) que la droite 
verticale $x_1=\cdots =x_n=0$ est incluse dans $D$ et coupe $V$ en un nombre fini de points n'appartenant 
pas à l'hyperplan à l'infini $z=0$. Ceci reste valable pour les droites voisines de $x=0$ 
ce qui permet d'utiliser les coordonnées affines $(x,y)=(x_1,..,x_n,y)$ au voisinage de $V$ et 
d'utiliser les coordonnées $(a,b)=(a_1,..,a_n,b_1,..,b_n)$ pour paramétrer l'ouvert affine de 
$(\pp^{n+1}(\cp))^*$ constitué des droites $L_{(a,b)}$ voisines de la droite verticale $x=0$, o\`u 
$L_{(a,b)}$ est la droite projective d\'efinie en coordonn\'ees affines $(x,y)$ par les 
\'equations 
$$
L_i(x,y,a,b)=x_i-a_i y-b_i=0\,,\ i=1,...,n.
$$
Dans ce cas, en supposant que $V=\{f(x,y)=0\}$ , on obtient 
(en utilisant les notations standard, voir par 
exemple \cite{pa:gnus}, pour d\'esigner les courants résiduels de Coleff-Herrera 
dans le cas intersection complète) l'expression locale de $T$~:
$$
T=\Phi\land df\wedge \Big( 
\bigwedge\limits_{i=1}^n d_{(x,y,a,b)}L_i\Big) 
\wedge \overline\partial \Big(\frac{1}{f}\Big)\wedge 
\Big(\bigwedge\limits_{i=1}^n 
\bar{\partial}_{(x,y,a,b)} \Big(\frac{1}{L_i}\Big)
$$
\vskip 2mm
\noindent
Examinons deux cas particuliers. 
\vskip 2mm
\noindent
{\bf 1.} {\it Le cas o\`u $\Phi$ est une fonction sur $V$} 
\vskip 2mm
\noindent
Si la forme $\Phi$ est une $0$-forme, c'est-\`a-dire 
une fonction m\'eromorphe sur $V$ ayant au voisinage de $V$ 
dans $D$ une expression globale $\Phi= h/g$ où $h$ et $g$ sont holomorphes en les coordonn\'ees affines 
$(x,y)$ et le lieu des zéros de $g$ coupe proprement $V$, l'expression de la trace de $\Phi$ au voisinage de $t=0$ dans $D^*$ est, 
en terme de calcul r\'esiduel~: 
$$
{\rm Tr}_V \Big(\frac{h}{g}\Big)
=\Bigg\langle 
\Big[\frac{1}{g}\Big]\bar{\partial}\Big(\frac{1}{f}\Big)\land
\bigwedge_{i=1}^n 
\overline {\partial}_{(x,y)}\Big(\frac{1}{L_i}\Big)\,,\, 
h\, J(f,L)dx\land dy \Bigg\rangle
$$
où 
$$
J(f,L)=\sum_1^n a_i\partial_{x_i}f+\partial_yf 
$$
est le jacobien de l'application $(x,y)\rightarrow (f,L_1,..,L_n)(x,y)$. 
Puisque les fonctions $(f,L_1,..,L_n)$ définissent une intersection complète 
et ont un nombre fini de zéros communs (dépendant holomorphiquement des paramètres $(a,b)$) 
et que la fonction $h$ est holomorphe, 
ce calcul de résidu global coïncide avec celui du résidu global de Grothendieck et on notera plus symboliquement
$$
{\rm Tr}_V \,\Big(\frac{h}{g}\Big)=
{\rm Res}\, \left[
\begin{matrix} 
{\frac{h}{g}}\, J(f,L) dx\land dy \cr 
f,L_1,..,L_n 
\end{matrix}
\right]\,. \eqno (\dag) 
$$
\vskip 2mm
\noindent
{\bf 2.} {\it Le cas o\`u $\Phi$ est une forme méromorphe de degré maximal sur $V$}
\vskip 2mm
\noindent 
Puisque $V$ ne contient pas la 
droite verticale $x=0$, la fonction $\partial_y f$ n'est pas identiquement nulle sur $V$ et l'on peut 
se ramener au cas o\`u $\Phi$ se pr\'esente sous la forme $\Phi=m(x,y)dx$, $m$ \'etant une fonction 
méromorphe sur $V$. Dans ce cas, on fait agir le courant $T$ sur des formes-test de $D\times D^*$ 
de bidegré $(n,2n)$ (en les variables $(a,b)$ ici) et on a alors l'expression suivante pour la forme trace 
${\rm Tr}_V\, (\Phi)$~: 
$$
{\rm Tr}_V \, (mdx)
=\sum_{k=0}^n {\rm Res}\, 
\left[\begin{matrix} 
m\, y^k\, \partial_yf\, dx\land dy \cr   
f,L_1,..,L_n
\end{matrix}
\right]
\Big(\sum_{|I|=k\,,\, |J|=n-k\,,\, {I\cap J =\emptyset}} \, \pm da_I\land db_{J}\Big)
\eqno(\dag\dag) 
$$
o\`u $I$ et $J$ d\'esignent des multi-indices ordonn\'es de $\{1,...,n\}$, $|I|$ et $|J|$ leurs 
cardinaux, et 
$$
da_I \wedge db_J := \bigwedge\limits_{l=1}^k da_{i_l} \wedge \bigwedge_{l'=1}^{n-k} db_{j_{l'}} 
$$
si $I=\{i_1,...,i_k\}$ et $J=\{j_1,...,j_{n-k}\}$. 

\begin{rem}
{\rm On retrouve ici une formule prouvée par P.A. Griffiths dans \cite{G:gnus} 
mais l'écriture $(\dag\dag)$ grâce au formalisme du calcul r\'esiduel 
met en \'evidence le fait que l'on n'ait pas \`a se soucier du fait que $V$ 
soit lisse ou non au voisinage des points d'intersection avec les droites 
correspondant aux points de $D^*$.
}
\end{rem}
\vskip 2mm
\begin{rem} 
{\rm A priori, $V$ n'est pas donné comme le lieu d'annulation d'une seule fonction. 
Par contre, si $V$ est un germe d'ensemble analytique en un point $(0,y_0)$ de 
la droite verticale distinct du point \`a l'infini, on a alors, en coordonn\'ees affines 
$(x,y)$ dans $D$ au voisinage de $(0,y_0)$,   
$V=\{f=0\}$ pour un germe de fonction holomorphe $f$ et 
les formules $(\dag)$ ou $(\dag\dag)$ s'appliquent à condition d'introduire au numérateur 
des symboles r\'esiduels une fonction plateau valant identiquement $1$ dans un voisinage 
suffisament petit de $(0,y_0)$. La trace est alors un germe de forme méromorphe en $(a,b)=(0,0)$.
}

\end{rem}
\vskip 2mm
\noindent
Un outil majeur du calcul r\'esiduel est le th\'eor\`eme de dualit\'e~: 
si $D\subset\cp^{n+1}$ est un domaine d'holomorphie, $(f_1,...,f_{k})$ 
une intersection complète dans $D$ et $g$ une fonction holomorphe dans $D$ telle 
que codim$\{g=0\}\cap \{f_1=\dots =f_k=0\}>k$, une fonction $h$ holomorphe 
dans $D$ est dans l'idéal engendré par les $f_i$ si et seulement si le 
courant r\'esiduel  
$$
h\Big[\frac{1}{g}\Big]\bar{\partial}(\frac{1}{f_1})\land...\land\bar{\partial}(\frac{1}{f_{k}})
$$
est nul sur l'espace des formes test de type $(n+1-k,n+1-k)$ $\overline\partial$-ferm\'ees au voisinage de 
$\{f_1=\dots=f_k=0\}$. Ce théorème s'applique notamment 
dans les anneaux locaux des germes de fonctions holomorphes. Le pendant alg\'ebrique de 
ce r\'esultat (dans le cadre \'el\'ementaire du calcul r\'esiduel en une variables) 
est une cons\'equence imm\'ediate de l'algorithme de division euclidienne~: 
si $F\in \cp [Y]$, un polyn\^ome $H$ de $\cp [Y]$ est divisible par $F$ si et seulement si 
$$
{\rm Res}\, \left[\begin{matrix}
Y^k\, H\, dY \cr F 
\end{matrix}
\right]=0\,,\qquad \forall k=0,..., \deg F-1\,. 
$$
 
\vskip 4mm
\noindent

\section{Construction du couple $(F,H)$ et nouvelle écri\-ture de la trace}

Nous allons dans ce paragraphe \'etablir 
une formule analogue à $(\dag\dag)$ dans le cas où $V$ est une 
réunion finie de germes d'hypersurfaces analytiques irr\'eductibles coupés (mais non inclus) 
par une droite que l'on supposera être la droite verticale. Cette formule fera encore intervenir des 
r\'esidus globaux de fractions rationnelles en une variable et l'on esquivera ainsi le 
recours aux fonctions plateau, ce qui fournira une expression diff\'erente de la trace ${\rm Tr}_V[mdx]$, 
d\'ependant du {\it degr\'e vertical} de $V$ (c'est-\`a-dire du nombre de points 
d'intersection de $V$ avec une droite g\'er\'erique $L_{(a,b)}$ pour $(a,b)$ voisin de $(0,0)$). 
On montrera pour cela l'équivalence, \`a degr\'e vertical $d$ pr\'ecis\'e, 
de la donn\'ee du couple $(V,\Phi)$ avec un couple de polynômes $(F,H)$ à une variable 
de degr\'es respectifs $d$ et $d-1$ et à coefficients dans le corps des germes de fonctions 
m\'eromorphes en $(a,b)$ \`a l'origine de $\cp^{2n}$.  
\vskip 2mm
\noindent
On notera ${\cal O}$ l'anneau factoriel 
(resp. ${\cal M}$ le corps) des germes de fonctions holomorphes 
(resp. méromorphes) à l'origine $(a,b)=(0,0)$ de $\cp^{2n}$. 
\vskip 2mm
\noindent
L'aboutissement de ce paragraphe est de caractériser ``le plus algébriquement possible'' 
à quelle condition un germe de forme méromorphe en un point de la grassmannienne $G(n+1,1)$ 
(correspondant par exemple \`a la droite $x=0$)   
est la trace d'un germe de $n$-forme méromorphe 
sur un "cycle de germes d'hypersurfaces analytiques" 
de "degré vertical ensembliste" donné (on expliquera plus loin 
ce vocabulaire). Donnons tout d'abord deux d\'efinitions dans un contexte plus g\'en\'eral 
que celui des sous-vari\'et\'es lin\'eaires, ce en vue de g\'en\'eralisations ult\'erieures de ce 
travail. 
\vskip 2mm
\noindent
\begin{deff} 
Soit $A$ une sous-variété algébrique de $\pp^{n+1}(\cp)$ de codimension pure $p$. 
On appelle germe de cycle analytique effectif intersectant proprement $A$ toute combinaison formelle à coefficients entiers 
positifs et finie de germes $\gamma$ d'espaces analytiques irréductibles de dimension $p$ en des points de $A$ 
avec $\dim ({\rm supp}\, (\gamma)\cap A)=0$. On note $\mathcal{V}_A$ l'ensemble des germes de cycles analytiques intersectant 
proprement $A$. Tout élément $V$ de  $\mathcal{V}_A$ peut se décomposer sous la forme 
$V=\sum_{P\in A}V_P$ où $V_P=\sum_i k_iV_{P,i}$ et les $V_{P,i}$ sont 
des germes d'ensembles analytiques irréductibles en $P\in A$ (en nombre fini), les points $P$
\`a consid\'erer \'etant aussi en nombre fini.  
On note $|V|$ l'ensemble analytique $\bigcup\limits_{P,i}V_{P,i}$ 
que l'on appelle le support de $V$.
\end{deff}
\vskip 2mm
\noindent
\begin{deff}
Si $V=\sum_{P\in A}V_P$ et $q\in \{0,...,n\}$, on note $M^q(|V|)$ l'ensemble des $(q,0)$-formes 
méromorphes sur $|V|$. Se donner une telle forme est équivalent à se donner 
un germe de forme méromorphe $\Phi_P$ sur chaque $|V_P|$, c'est-\`a-dire un germe de 
forme m\'eromorphe dans l'espace $\pp^{n+1}(\cp)$ au voisinage de chaque point $P$ de $A\cap|V|$.  
\end{deff}
\vskip 2mm
\noindent
On s'intéresse ici au cas où $A$ est la droite verticale $\Delta=L_{(0,0)}~: x=0$. 
On notera $\mathcal{V}$ le sous-ensemble de $\mathcal{V}_{\Delta}$ 
constitué des combinaisons de germes d'hypersurface en des points $P=(0,y_P)$ de $\Delta$
distincts du point \`a l'infini de cette droite.  
\vskip 2mm
\noindent
Soit $V$ un germe d'hypersurface analytique irréductible en un point $P=(0,y_P)\in\Delta$ intersectant proprement 
$\Delta$ en ce point. Dans ce cas, $V=\{f=0\}$  où $f\in\cp\{x,y-y_P\}$ est un germe de fonction holomorphe 
réduit en $P$ (vérifiant donc $({\rm rad}\, (f))=(f)$). Par hypothèse, la fonction $f(0,y)$ n'est pas identiquement nulle 
et on note $d$ son ordre d'annulation en $y=y_P$ . Par continuité, $d$ est aussi le nombre de zéros de 
la fonction holomorphe $y\rightarrow f(x,y)$ pour $x$ voisin de zéro. On appellera $d=\deg(V)$ 
le {\it degré vertical} de $V$. 
Le degré vertical d'un cycle $V=k_1V_1+..+k_sV_s$ sera alors par définition 
$\deg(V)=k_1 \deg(V_1)+..+k_s \deg(V_s)$. L'ensemble $\mathcal{V}$ a naturellement une structure de 
semi-groupe gradué avec la graduation par le degré vertical (on rajoute l'ensemble vide de degré $0$ pour avoir un élément neutre). 
Tout élément admet une décomposition unique en somme d'éléments irréductibles.
\vskip 2mm
\noindent
\begin{rem}
{\rm Si $V$ est un germe en un point $P\in\Delta$, l'entier $d$ est le degré du revêtement analytique 
$(V,\pi,I)$ o\`u $I$ est un ouvert suffisamment petit de $\cp^n$ 
centré en  $x=0$ et $\pi$ est la projection verticale $(x,y)\rightarrow x$. 
Notons que $d$ n'est pas forcément le degré global de $V$ car il peut varier d'une projection \`a l'autre, 
la projection choisie ne correspondant pas n\'ecessairement \`a une projection g\'en\'erique~: par exemple le 
degré vertical du germe en $0$ de l'ensemble $V=\{y^2-x^3=0\}$ est $2$ alors que le degré de $V$ est 
$3$ pour une projection g\'en\'erique. 
}
\end{rem}
\vskip 2mm
\noindent 
Quand on parle d'une équation $\{f=0\}$ d'un germe $V$, on suppose que l'écriture de $f$ prend en compte 
les multiplicités de chacune des branches de $V$.
\vskip 2mm
\noindent
\begin{rem} 
{\rm La définition de la trace avec multiplicité ({\it i.e.} la fonction $f$ n'est pas forcément réduite) 
gr\^ace aux formules $(\dag)$ ou $(\dag \dag)$ revient à la définition usuelle~:
$$
{\rm Tr}_{k_1V_1+...+k_sV_s}\, (\Phi)={\rm Tr}_{V_1+\cdots +V_s}\, (\Phi')
$$
où $\Phi'$ est la forme méromorphe sur $|V|$ valant $k_i\Phi$ sur $V_i$.
C'est seulement dans la construction de $F$ (rendant compte de $V$) que l'on fera entrer en jeu les multiplicit\'es 
impliqu\'ees dans la d\'efinition de $V$. On se restreindra ensuite aux ensembles 
${M}(|V|)$ (resp. ${M}^n(|V|)$) des fonctions (resp. des $n$-formes) 
m\'eromorphes sur $|V|$ pour construire $H$.
}
\end{rem}
\vskip 2mm
\noindent
On se ramène dans un premier temps au calcul résiduel alg\'ebrique en une variable grâce au théorème de préparation 
de Weierstrass. Soit $V=\{f=0\} \in\mathcal{V}$ un germe irréductible en un point $P\in\Delta$, de degr\'e vertical $d$. 
Soit $\wt{f}(y,a,b)=f(ay+b,y)\in {\cal O}\, \{y-y_P\}$. La fonction $\wt{f}$ est irréductible et 
régulière de degré $d$ en $y=y_P$. Par le théorème de préparation de Weierstrass, 
il existe un unique polynôme irréductible $Q\in {\cal O}\, [Y]$ de degré $d$ et une fonction analytique $u\in {\cal O}\, \{y-y_p\}$ 
inversible tels que~:
$$
\wt{f}(y,a,b)=Q(y,a,b)u(y,a,b)
$$
avec 
$Q(y,0,0)=(y-y_P)^d$ et $u(y_P,0,0)\ne 0$. On a alors le lemme suivant, permettant  
pr\'ecis\'ement de nous ramener au calcul r\'esiduel en une variable~: 
\vskip 2mm
\noindent
\begin{lem} 
Soit $V=\{f=0\} \in\mathcal{V}$ un germe irréductible de degr\'e vertical $d$ et 
$Q$ le polynôme irréductible de 
${\cal O}\, [Y]$ défini ci-dessus. 
Soit $h\in {M}(|V|)$ un germe de fonction méromorphe sur $|V|$ (c'est-\`a-dire un germe de 
fonction m\'eromorphe au point $(0,y_P)$ de $\pp^{n+1}(\cp)$)~; on a alors
$$
{\rm Tr}_V \, (h) = {\rm Res}\, \left[\begin{matrix} h(ay+b,y)\, \partial_yQ(y,a,b)\, dy\cr   Q(y,a,b) \end{matrix}\right].
$$
\end{lem}
\vskip 2mm
\noindent {\bf Preuve.} 
On introduit une fonction plateau $\theta$ valant $1$ au voisinage de $V$. Ainsi, par $(\dag)$ et 
au vu de la remarque $2$, l'expression de la trace devient, si $J(f,L):=J(f,L_1,\cdots,L_n)$, 
$$
{\rm Tr}_V (h)={\rm Res}\, 
\left[\begin{matrix}\theta(x,y)\, h(x,y)\, J(f,L)\, dy\land dx_1\land ..\land dx_n\cr f,x_1-a_1y-b_1,\dots ,x_n-a_ny-b_n
\end{matrix}\right]
$$
On remarque que la fonction de $(x,y,a,b)$
$$
(x,y,a,b)\to \theta \, h\, J(f,L)\, (x,y,a,b)-\theta \, h\, J(f,L)(ay+b,y,a,b)
$$
s'écrit comme combinaison linéaire des $L_1,..,L_n$ à coefficients semi-méro\-morphes en $(0,y_P,0,0)$, 
le lieu polaire de ces coefficients ne contenant pas la droite $x-ay-b=0$ pour $a,b$ g\'en\'eriques au voisinage de 
$(0,0)$. En dehors du lieu polaire de $h$, le théorème de dualité s'applique~; ainsi sur l'ouvert constitué des 
$(a,b)$ pour lesquels la droite $L(a,b)$ ne passe pas par l'ensemble fini des points communs à $V$ et au 
lieu polaire de $h$, on a l'égalité
\begin{eqnarray*}
{\rm Tr}_V \, (h) &=& {\rm Res}\, \left[\begin{matrix}\theta\, h(x,y)\, J(f,L)(x,y,a,b)\, dy\land dx \cr  
f,L_1,\dots ,L_n\end{matrix}\right] \\
&=& {\rm Res}\, \left[\begin{matrix}\theta (ay+b,y)\, h(ay+b,y)\, 
J(f,L)(ay+b,y,a,b)\, dy\land dx\cr  f(ay+b,y),L_1,..,L_n \end{matrix}\right]\,, 
\end{eqnarray*}
d'o\`u l'on d\'eduit l'égalité dans ${\cal M}$. Cette expression nous permet de nous lib\'erer des variables $x$ 
et d'exprimer la trace comme un calcul de r\'esidu de forme m\'eromorphe en une variable $y$, 
\`a savoir  
$$
{\rm Tr}_V \, (h)= {\rm Res}\, 
\left[\begin{matrix} \theta (ay+b,y)\, h(ay+b,y)\, J(f,L_1,..,L_n)(ay+b,y,a,b)\, dy\cr f(ay+b,y)
\end{matrix}
\right]
$$
Un rapide calcul montre que l'on a l'égalité suivante~:
$$
J(f,L_1,..,L_n)(ay+b,y,a,b)=\Big(\sum_1^n a_i\partial_{x_i}f+\partial_yf\Big)(y,a,b)=\partial_y\wt{f}\, (y,a,b)~; 
$$
puisque $\wt{f}(y,a,b)=Q(y,a,b)u(y,a,b)$, on a
\begin{eqnarray*} 
{\rm Tr}_V \, (h) &=& {\rm Res}\, \left[\begin{matrix} 
\theta(ay+b,y)\, h(ay+b,y)\, \partial_y\wt{f}(y,a,b)dy\cr \wt{f}(y,a,b) \end{matrix}\right] \\
&=&  {\rm Res}\, 
\left[\begin{matrix} \theta (ay+b,y)\, h(ay+b,y)\, \partial_yQ\, dy\cr Q(y)\end{matrix}\right]+ \\
& &   {\rm Res}\, \left[\begin{matrix} 
\theta (ay+b,y)\, h(ay+b,y)\, \partial_y u\, dy\cr  u \end{matrix}\right]~; 
\end{eqnarray*}
la fonction $u$ étant inversible sur le support de $y\to p(ay+b,y)$, 
le deuxième terme de la somme de droite est nul~; de plus les zéros de $Q$ étant localisés au voisinage de l'origine, 
la fonction plateau devient inutile dans le premier terme de la somme de droite, ce qui 
fournit la formule souhait\'ee pour la trace ${\rm Tr}_V\, (h)$. \hfill $\square$ 
\vskip 2mm
\noindent 
Si maintenant $V=V_1+\cdots + V_k\in \mathcal{V}$ est une somme de germes irréductibles 
$V_i=\{f_i(x,y)=0\}$ distincts deux à deux, on définit les fonctions $Q_i\in {\cal O}[Y]$ 
associées aux $f_i$ comme précédemment et, en posant $Q=Q_1\cdots Q_k$, on obtient~:
$$
{\rm Tr}_V \, (y^k)
=\sum_i {\rm Res}\, \left[\begin{matrix} y^k\, \partial_y\, Q_i \, dy\cr  Q_i \end{matrix}\right]=
{\rm Res}\, \left[\begin{matrix} y^k\, \partial_yQ \, dy\cr Q \end{matrix} \right]\,. 
$$
On notera $u_k$ ces fonctions de $(a,b)$ qui sont des germes de fonctions holomorphes par le théorème d'Abel. 
On remarque que $u_0$ est le degré de $Q$ qui est aussi le degré vertical $d$ de $V$. 
On a alors le deuxième lemme crucial dans l'exposé~:
\vskip 2mm
\noindent
\begin{lem} La matrice $d\times d$ suivante $U$ \`a coefficients dans ${\cal M}$
$$
A=\left(\begin{matrix}
u_{0} & u_{1} & \cdots & u_{d-1} \cr
u_1 & u_2 & \cdots & u_d \cr 
\vdots & \vdots & \vdots & \vdots \cr
u_{d-1} & u_{d} & \cdots & u_{2d-2}
\end{matrix}
\right)  
$$ 
est non dégénérée sur le corps ${\cal M}$. Plus précisément, on a
$$
{\rm Det}\,A (a,b)\,=\, {\rm Disc}\, Q (a,b)
$$
où ${\rm Disc}\,Q$ est le discriminant de $Q$.
\end{lem}
\vskip 2mm
\noindent 
{\bf Preuve.} Dans un premier temps, il suffit de montrer que le noyau de l'endomorphisme $U$ associé à cette matrice est nul. Soit 
$$
\sigma=(\sigma_{0},...,\sigma_{d-1})\in \mathcal{M}^d 
$$
un $d$-uplet d'\'ements de ${\cal M}$ tel que $U(\sigma)=0$. 
En terme de calcul résiduel, cette condition se traduit par
$$
{\rm Res}\, \left[\begin{matrix} Y^k(\sigma_{d-1}Y^{d-1} + \cdots + \sigma_1 Y+ \sigma_0)\, \partial_yQ (Y,a,b)\, dY
\cr Q(Y,a,b) \end{matrix} \right]=0\,,\ k=0,\dots,d-1 
$$ 
ce qui par le théorème de dualité se traduit par l'appartenance du 
polynôme $(\sigma_{d-1}Y^{d-1}+\cdots + \sigma_1 Y +\sigma_0)\partial_yQ$ à l'idéal engendré par $Q$ dans 
${\cal M}[Y]$. Puisque $Q$ est sans facteurs multiples et pour une raison de degré, 
on voit par le lemme de Gauss que l'unique solution est $\sigma=(0,...0)$, ce qui prouve la première partie du lemme. On voit par ce biais que ${\rm Det}\, A(a_0,b_0)=0$ est équivalent au fait que les polynômes $\partial_yQ(Y,a_0,b_0)$ et $Q(Y,a_0,b_0)$ ont une racine commune (pour $(a,b)=(a_0,b_0)$ fixés). En fait, pour $(a,b)$ générique, le polynôme $Q$ a $d$ racines $y_1(a,b),...,y_d(a,b)$ distinctes. Il suffit alors de remarquer que $A=S^tS$ où $S=(y_i^j)_{0\le i,j\le d-1}$ est la matrice de Vandermonde des racines de $Q$ ce qui montre la deuxième partie du lemme.
\hfill $\square$
\vskip 2mm
\noindent
\begin{rem}
{\rm On a privilégié la démonstration du lemme par le théorème de dualité qui est une approche plus globale (on ne travaille pas génériquement) qui s'adapte pour des situations plus compliquées. Quelle que soit l'approche, on comprend bien que les droites $L(a,b)$ pour lesquelles ${\rm Det} A(a,b)=0$ sont exactement les droites pour lesquelles $f(ay+b,y)$ a une racine double: ce sont les ``mauvaises'' droites, c'est à dire celles qui passent par le lieu singulier de $V$ ou qui sont tangentes à $V$. }
\end{rem}
\vskip 2mm
\noindent
On note ${\cal U}\, [Y]$ le sous-ensemble de ${\cal O}[Y]$ constitué  
des polynômes unitaires vérifiant les $n$ relations
$$
\partial_{a_i} F -Y\partial_{b_i} F \in (F)\,,\qquad i=1,\cdots,n  
\eqno(*)
$$
\vskip 2mm
\noindent
\begin{lem}
L'ensemble ${\mathcal U}\, [Y]$ admet une structure de semi-groupe (multiplicatif) 
gradué avec la graduation naturelle qui de plus hérite de la factorialité de l'anneau ${\cal O}[Y]$.
\end{lem}
\vskip 2mm
\noindent
{\bf Preuve.}  Si $F_1$ et $F_2$ sont deux facteurs de $F$ vérifiant 
$$
\partial_{a_i} F_j -Y\partial_{b_i} F_j \in (F_j)\quad j=1,2\,, 
$$
on a alors 
\begin{eqnarray*}
\partial_{a_i}  (F_1F_2) -Y\partial_{b_i} (F_1F_2) &=& F_2(\partial_{a_i} F_1 -Y\partial_{b_i} F_1)+F_1(\partial_{a_i} F_2 -Y\partial_{b_i} F_2)\\
 &\in& (F_1F_2)\,, 
\end{eqnarray*}
et ce pour tout $i=1,..,n$. La propriété $(*)$ est donc compatible avec la multiplication, 
ce qui donne à $\mathcal{U}\,[Y]$ une structure de semi-groupe (multiplicatif) 
gradué avec la graduation naturelle. Il reste à montrer la factorialité. 
Soit $F\in \mathcal{U}\,[Y]$ et soit $F=F_1^{k_1}\cdots F_s^{k_s}$ sa décomposition dans l'anneau factoriel 
${\cal O}[Y]$. Il suffit de montrer que tout facteur irréductible 
$F_i$ de $F$ vérifie $(*)$~; on le montre pour $F_1$~;  
le polynôme $F$ se factorise en $F=F_1^{k_1} P$ où $k_1\in \mathbb{N}^*$ et 
$P$ est premier avec $F_1$ et l'on a 
\begin{eqnarray*} 
\partial_{a_i} F -Y\partial_{b_i} F &=& F_1^{k_1}(\partial_{a_i} P -Y\partial_{b_i} P)
+kF_1^{k_1-1}P(\partial_{a_i} F_1 -Y\partial_{b_i} F_1)\\
&\in & (F_1^{k_1} P)~; 
\end{eqnarray*} 
ceci implique $k_1F_1^{k_1-1}P(\partial_{a_i} F_1 -Y\partial_{b_i} F_1)\in (F_1^{k_1})$ et ce pour tout $i$. 
Or, $P$ est premier avec $F_1$ par hypothèse et on a donc 
$\partial_{a_i} F_1 -Y\partial_{b_i} F_1\in (F_1)$, ce qui ach\`eve la preuve du lemme $3$. 
\hfill $\square$ 
\vskip 2mm
\noindent
\begin{prop} 
Il existe un isomorphisme $\Pi$ de semi-groupes gradués entre les ensembles $\mathcal{V}$ et $\mathcal{U}\, [Y]$~; 
les branches irréductibles de $V$, leur degré vertical et leur multiplicité sont 
respectivement en correspondance avec les facteurs irréductibles de $F$, leur degré et leur multiplicité 
dans la décomposition de $F$ dans l'anneau factoriel ${\cal O}[Y]$. 
\end{prop} 
\vskip 2mm
\noindent
{\bf Preuve.} Il suffit de construire l'image d'un seul germe irréductible et 
l'image réciproque d'un polynôme irréductible, puis de montrer que $\Pi$ est un homomorphisme de semi-groupes gradués. 
Soit donc $V$ un tel \'el\'ement de $\mathcal{V}$ de degr\'e vertical $d$. L'image de $V$ par $\Pi$ sera par définition
$$
F:=\Pi(V)=Y^d-\sigma_{d-1}Y^{d-1}+...+(-1)^{d-1}\sigma_0\,, 
$$
où $(\sigma_{d-1},...,\sigma_{0})\in {\cal M}^d$ est l'unique solution 
du système linéaire $(S)$ de $d$ équations à $d$ inconnues 
\begin{eqnarray*} 
\begin{matrix}
u_{d-1}\sigma_{d-1} &+& \cdots &+&(-1)^{d-1}u_{0}\sigma_{0} = u_d \\
&\vdots &\vdots &\vdots & \vdots \\
u_{2d-2}\sigma_{d-1} &+&\cdots &+&(-1)^{d-1}u_{d-1}\sigma_{0} = u_{2d-1} 
\end{matrix} 
\end{eqnarray*} 
qui est un syst\`eme de Cramer d'apr\`es le lemme $2$. Le fait que $\sigma$ soit solution 
de ce syst\`eme de Cramer est équivalent \`a ce que le polyn\^ome $F$ construit ci-dessus \`a partir de $\sigma$ v\'erifie 
$$
{\rm Res}\, \left[ \begin{matrix} y^k\, F\partial_yQdy\cr Q \end{matrix}\right] =0\,, \quad k=0,\cdots,d-1\,. 
$$
Par le théorème de dualité on a donc $F\, \partial_yQ\in(Q)$. Ici, $Q$ est irréductible et unitaire par hypothèse, 
on a donc $Q=F$ pour des raisons de degré. L'unique solution du syst\`eme de Cramer $(S)$ fournit 
donc le polyn\^ome de Weierstrass $Q \in {\cal O}[Y]$ associ\'e \`a la fonction $\wt{f}$.  
Par définition même de $\wt{f}$, on a pour tout $i=1,\dots,n$, 
$$
0=\partial_{a_i} \wt{f}-Y\partial_{b_i}\wt{f}=u(\partial_{a_i} F -Y\partial_{b_i} F)+F(\partial_{a_i} u -Y\partial_{b_i} u)
$$
et $F$ divise $\partial_{a_i} F -Y\partial_{b_i} F$ dans 
${\cal O}[Y]$ pour tout $i=1,..,n$, ce qui montre que le polynôme $F$ appartient \`a $\mathcal{U}\, [Y]$.  
De plus, par construction m\^eme de $Q$, $\Pi$ est un homomorphisme de semi-groupes gradués entre 
$\mathcal {V}$ et $\mathcal {U}\, [Y]$.
\vskip 2mm
\noindent
Prouvons maintenant la {\it surjectivité} de $\Pi$. 
Il suffit de montrer que tout polynôme irréductible de $\mathcal {U}\, [Y]$ 
s'\'ecrit $\Pi(V)$ o\`u $V\in \mathcal{V}$. Soit $F\in \mathcal{U}\, [Y]$ 
un polyn\^ome irréductible de degré $d$. Du fait de l'irr\'eductibilit\'e de 
$F$, le polyn\^ome $F(Y,0,0)$ ne peut avoir qu'une seule racine $y_P$ et 
s'ecrit donc $F(Y,0,0)= (Y-y_P)^d$. On considère alors la fonction 
$(x,y,a)\to  G(x,y,a):=F(y,a,x-ay)$ qui, vue comme un élément de 
$\cp\{x,y-y_P,a\}$, est régulière de degré $d$ en $y=y_P$. 
Puisque $F$ est dans $\mathcal {U}\, [Y]$, donc v\'erifie $(*)$ et est unitaire, 
la fonction $G$ appartient \`a l'id\'eal engendr\'e par $\partial_{a_i} G$ 
dans $\cp \{x,y,a\}$~; or l'ordre d'annulation en $a_j=0$ 
de la fonction $a_j\rightarrow G(x,y,a)$ est donné par l'intégrale 
$$
\frac{1}{2i\pi}\int_{|a_j|=\ep}\frac{\partial_{a_j} G \, (x,y,a)}{G(x,y,a)}\, da_j
$$
(pour $\ep$ assez petit) et vaut zéro pour tout $(x,y)$ voisin de $(0,y_P)$ 
puisque la fonction sous l'intégrale est holomorphe pour $\ep$ suffisament petit. 
Par utilisations successives du théorème de préparation de Weierstrass (on \'elimine les 
variables $a_j$, $j=1,...,n$ les unes apr\`es les autres), on montre 
alors l'existence d'un germe de fonction holomorphe unique (à inversible près) $f\in\cp\{x,y-y_P\}$ 
s'annulant en $(0,y_P)$ et d'un germe de fonction holomorphe inversible $q\in\cp\{x,y-y_P,a\}$ tels que 
$$
G(x,y,a)=f(x,y)q(x,y,a)\,, 
$$
où $f$ est régulière de degré $d$ au point $(0,y_P)$. On a alors 
$$
\wt{f}(y,a,b)=G(ay+b,y,a)q^{-1}(ay+b,y,a)=F(y,a,b)q^{-1}(ay+b,y,a)
$$
Il suffit de poser $V=\{f=0\}$ et $V$ est alors un germe d'ensemble analytique en 
$P=(0,y_P)\in\Delta\setminus H_\infty$ tel que $F=\Pi(V)$. Ce germe est irréductible sinon $F$ serait réductible. 
\vskip 2mm
\noindent
Prouvons maintenant l'{\it injectivit\'e} de $\Pi$. Grâce à la factorialité de $\mathcal{U}\, [Y]$, on se ram\`ene   
\`a montrer l'injectivité pour des germes irréductibles. 
Soient $V_1=\{f_1=0\}$ et $V_2=\{f_2=0\}$ deux germes analytiques irréductibles en deux points 
$P_1=(0,y_{P_1})$ et $P_2=(0,y_{P_2})$ de $\Delta\setminus H_\infty$ 
tels que $\Pi(V_1)=\Pi(V_2)=F$. Dans ce cas on a 
$\wt{f_1}(y,a,b)=F(y,a,b)u_1(y,a,b)$ dans ${\cal O}\, \{y-y_{P_1}\}$ et 
$\wt{f_2}(y,a,b)=F(y,a,b)u_2(y,a,b)$ dans ${\cal O}\, \{y-y_{P_2}\}$. 
Puisque $F$ est irréductible dans $\mathcal{U}\, [Y]$, $F(y,0,0)=(y-y_0)^d$, 
ce qui montre que $P_1=P_2$~; ainsi les deux fonctions  $\wt{f_1}$ et $\wt{f_2}$ 
sont définies au voisinage d'un même point et sont égales à un inversible près~; 
par conséquent $f_1(x,y)$ et $f_2(x,y)$ également. On a donc bien $V_1=V_2$. 
\vskip 2mm
\noindent
La proposition $1$ est ainsi d\'emontrée. \hfill $\square$ 
\vskip 2mm
\noindent
\begin{rem} 
{\rm On peut noter que si $\{f=0\}$ est l'\'equation de $V$ 
(les multiplicit\'es \'etant prises en compte), les coefficients de $F=\Pi(V)$ sont 
les fonctions symétriques élémentaires en les racines de $y\to \wt{f}(y,a,b)$ et 
peuvent s'exprimer grâce aux formules de Newton comme des polynômes 
en les polynômes symétriques de Newton des racines que 
sont les $u_i$ pour $i=0,..,d$. Il n'était donc pas nécessaire d'utiliser le 
système de Cramer $(S)$ pour démontrer la proposition 
mais l'utilité de ce type de construction appara{\^\i}tra par la suite.
}
\end{rem}
\vskip 2mm
\noindent
\begin{rem} 
{\rm Le fait que l'ensemble analytique $V$ contienne la droite verticale $\{x=0\}$ équivaut 
au fait que la fonction holomorphe $f$ soit dans l'idéal engendré par les $x_i$. Le polynôme $F$ 
alors obtenu serait dans l'idéal engendré par les polynômes  $Y+ b_i/a_i$, $i=1,...,n$~; ce  
polynôme vérifie toujours $(*)$ mais ses coefficients ne sont plus holomorphes et n'ont plus de 
signification en $(a,b)=(0,0)$~; ceci s'explique du fait que la droite $(a,b)=(0,0)$ est incluse dans $V$. M\^eme 
si ce cas s'av\`ere pathologique, on peut par contre montrer, d'une manière générale, 
que le lieu polaire (dans $(\cp^{2n},0)$) des coefficients d'un \'el\'ement de $\mathcal{M}\, [Y]$ v\'erifiant 
les conditions $(*)$ ne dépend en fait pas de $b$.}
\end{rem}
\vskip 2mm
\noindent 
On travaille maintenant ensemblistement, {\it i.e.} 
sur l'ensemble $\mathcal{V}_{\rm red}$ constitué des  
\'el\'ements $V=V_1+\cdots +V_k\in\mathcal{V}$, où les $V_i$ sont deux à deux distincts 
(et $F=\Pi(V)$ n'a pas de facteurs multiples).
\vskip 2mm
\noindent
\begin{prop} Soit $V\in \mathcal{V}_{\rm red}$ de degré vertical $d$ et $F=\Pi(V)$. 
Il existe une application bijective 
$$
\rho~: {\cal M} (V)\rightarrow \mathcal{M}_F\, [Y]
$$ 
entre l'ensemble $\mathcal{M} (V)$ des germes de fonctions méromorphes sur $V$ et 
le sous-ensemble $\mathcal{M}_F\, [Y]$ de $\mathcal{M}\, [Y]$ constitué des polynômes 
$H$ de degré $d-1$ \`a coefficients dans $\mathcal{M}$ vérifiant la relation
$$
\partial_{a_i} H -Y\partial_{b_i} H \in (F)\quad\forall\, i=1,..,n\,. \eqno (**)
$$
\end{prop}
\vskip 2mm
\noindent {\bf Preuve.} La preuve se fait en deux \'etapes~: nous traitons tout d'abord le 
cas d'un germe irr\'eductible $V$ (en un point $(0,y_P)$), puis le cas g\'en\'eral. 
\vskip 2mm
\noindent
{\it Le cas d'un germe irr\'eductible.} 
Soit $V$ un germe d'équation (réduite) $V=\{f=0\}$ de degré vertical $d$ et $F=\Pi(V)$. 
Soit $h$ un germe de fonction méromorphe sur $V$, $h\in \mathcal{M} (V)$. 
La construction de $\rho$ est analogue à celle de $\Pi$ mais on considère cette fois le système 
$(S_h)$ de $d$ équations à $d$ inconnues suivant (les $\tau_i$ sont les inconnues)~:
\begin{eqnarray*} 
\begin{matrix} u_{0}\tau_{0} &+& \cdots &+&u_{d-1}\tau_{d-1} = v_0 \\
&\vdots & \vdots &\vdots & \vdots \\
u_{d-1}\tau_{0} &+& \cdots &+&u_{2d-2}\tau_{d-1} = v_{d-1} 
\end{matrix} 
\end{eqnarray*} 
avec cette fois les $v_k$ (ce sont des \'el\'ements de $\mathcal{M}$ d'apr\`es le th\'eor\`eme d'Abel) 
d\'efinis par 
$$
v_k:={\rm Tr}_V\, (y^k\, h)={\rm Res}\, \left[\begin{matrix} 
y^k\, h(ay+b,y)\, \partial_yF\, dy\cr F(y,a,b)\end{matrix}\right] \,,\qquad k=0,...,d-1\,. 
$$
Soit $\tau:=(\tau_{0},..,\tau_{d-1})$ l'unique $d$-uplet solution du 
système (qui est de Cramer d'apr\`es le lemme $2$), qui cette fois est {\it a priori} 
un vecteur de germes de fonctions m\'eromorphes en $(a,b)$ \`a l'origine de $\cp^{2n}$. 
On pose alors 
$$
H(Y,a,b)=\tau_{d-1}(a,b)Y^{d-1}+\cdots + \tau_1(a,b) Y +\tau_0(a,b)\,. 
$$ 
Par construction même, on a~:
\begin{eqnarray*} 
{\rm Res}\, 
\left[
\begin{matrix} 
H(y,a,b)\, y^k\, \partial_yF(y,a,b)\, dy\cr
F(y,a,b)\end{matrix}\right]
&=&\sum_0^{d-1}\tau_i {\rm Res}\, \left[
\begin{matrix} y^{k+i}\, \partial_yF(y,a,b)\, dy\cr 
F(y,a,b) \end{matrix} \right]\\
&=&\sum_0^{d-1}\tau_iu_{k+i}=v_k
\end{eqnarray*} 
pour $k=0,...,d-1$. Ainsi, 
$$
{\rm Res}\, \left[\begin{matrix} 
y^k\, h(ay+b,y)\, \partial_yF\, dy\cr F(y,a,b)\end{matrix}\right] = {\rm Res}\, 
\left[
\begin{matrix} 
H(y,a,b)\, y^k\, \partial_yF(y,a,b)\, dy\cr
F(y,a,b)\end{matrix}\right]
$$
pour tout $k=0,...,d-1$. Puisque $d=\deg_Y F$, on en déduit que ceci est en fait vrai pour tout entier $k$ 
en effectuant la division euclidienne de $y^k$ par $F$ et en utilisant le théorème de dualité. 
Ainsi l'égalité reste vraie si on multiplie les deux numérateurs par un \'el\'ement 
quelconque de $\mathcal{O}\,\{y-y_P\}$. Or, comme $h$ est la restriction \`a 
$|V|$ d'un germe de forme m\'eromorphe au voisinage de $(0,y_P)$ dans $\pp^{n+1}(\cp)$, $h$ admet 
un d\'enominateur $(x,y)\to \xi(x,y)$ (d'ailleurs puissance d'un d\'enominateur universel $\delta$ ne d\'ependant que de 
$V$ d'apr\`es le th\'eor\`eme d'Oka) qui est un germe de fonction holomorphe au voisinage de $(0,y_P)$~;
si l'on pose $r(y,a,b)=\xi (ay+b,y)$, il r\'esulte des \'egalit\'es entre symboles r\'esiduels ci-dessus 
(pens\'ees dans $\mathcal{M}\, \{y-y_P\}$) et du th\'eor\`eme de dualit\'e que 
$(y,a,b) \to r(y,a,b) (H(y,a,b) - h(ay+b,y)) \partial_y F$ est dans l'id\'eal engendr\'e par $F$
dans $\mathcal{O}\, \{y-y_P\}$~; on en d\'eduit (car $F$ est irréductible par hypothèse) qu'en tant qu'\'el\'ement de $\mathcal{M}\, \{y-y_P\}$, $H$ 
s'\'ecrit 
$$
(y,a,b) \to H(y,a,b)= h(ay+b,y)+ \frac {q(y,a,b)}{r(y,a,b)} F(y,a,b) \eqno(\dag^3)  
$$
et puisque les fonctions $(\partial_{a_i}- y \partial_{b_i})(r)$ et      
$(\partial_{a_i}- y \partial_{b_i})(h(ay+b,y))$ pour $i=1,...,n$ sont identiquement nulles et que $F$ vérifie $(*)$, 
le polynôme $\partial_{a_i} H -y \partial_{b_i} H$ est donc dans l'id\'eal engendr\'e par $F$ dans $\mathcal{O}\, \{y-y_P\}$~; 
il en r\'esulte que $H$ v\'erifie le jeu d'\'equations $(**)$ et est donc dans $\mathcal{M}_F\, [Y]$. 
Cette construction nous permet donc d'associer \`a $h$ un \'el\'ement $\rho(h)$ de $\mathcal{M}_F\, [Y]$. 
\vskip 2mm
\noindent
L'injectivit\'e de $\rho$ r\'esulte de la remarque \'evidente suivante~: d'apr\`es $(\dag^3)$, la restriction \`a 
$V$ de $(x,y)\to H(y,a,x-ay)$ (consid\'er\'e comme germe 
d'une fonction de $x,y$ au voisinage de $(0,y_P)$) lorsque $a$ est fix\'e 
g\'en\'eriquement et voisin de $0$ est \'egale \`a $h$. 
\vskip 2mm
\noindent
Prouvons maintenant la surjectivit\'e de $\rho$. 
Soit $H\in\mathcal{M}_F\, [Y]$. Puisque $H$ vérifie les conditions $(**)$ et que $F$ est unitaire,
on a 
$$
\partial_{a_i} H-Y\partial_{b_i} H =q_i(a,b)F(Y,a,b)
$$
avec $q_i= -\partial_{b_i}\, \gamma_{d-1}$, $i=1,...,n$~; ainsi    
\begin{eqnarray*}
\partial_{a_i} [H(y,a,x-ay)] &=& q_i(a,x-ay)F(y,a,x-ay) \\
&=& q_i(a,x-ay)\, u(a,x,y)f(x,y)
\end{eqnarray*} 
où $u$ est inversible et $V=\{f=0\}$. Le lieu polaire de la fonction méromorphe 
$(x,y)\rightarrow q_i(a,x-ay)$ ne contient 
aucune branche de $V$, sinon $q_i(a,b)$ serait divisible (dans $\mathcal {M}\, [Y]$) 
par un facteur de $F(Y,a,b)$, ce qui est absurde. 
On peut donc restreindre à $V$ les deux membres de cette 
derni\`ere \'egalit\'e, ce qui donne~: 
$$
\partial_{a_i} [H(y,a,x-ay)]_{|V}=\partial_{a_i} [H(y,a,x-ay)_{|V}]=0\,, \quad i=1,...,n\,, 
$$
la première de ces égalités \'etant une conséquence du fait que l'équation de $V$ ne dépende pas de $a$. 
La fonction $H(y,a,x-ay)$ restreinte à $V$ ne dépend donc pas de $a$ et 
permet de définir un germe de fonction méromorphe $h$. 
Les coefficients de $H$ sont alors solutions du système de Cramer $(S_h)$ obtenu à partir de $h$, 
ce qui implique que $H=\rho(h)$.
\vskip 2mm
\noindent
{\it Le cas d'un cycle r\'eduit.} 
Soit maintenant  $V=V_{P_1}+ \cdots + V_{P_s}\in \mathcal {V}_{\rm red}$ un cycle r\'eduit 
de degré vertical $d$, où les $V_{P_j}$, $j=1,...,s$, 
sont des germes de cycles de dimension $n$ irr\'eductibles en les $s$ points $P_j$ distincts. 
Soit $F_j=\Pi(V_{P_j})$ le polynôme de degré $d_j$ associé à $V_{P_j}$ et $F=F_1\cdots F_s=\Pi(V)$ 
le polynôme de degré $d$ associé à $V$. On note $F_{[j]}=\prod_{l\ne j}F_l$, $j=1,...,s$~; 
le fait que les germes n'aient pas de branches communes 
implique que les $F_j$ sont premiers deux à deux~; par le théorème de B\'ezout,  
il existe $s$ polynômes $U_1,..,U_s\in \mathcal {\cal M} [Y]$ tels que $U_1F_{[1]}+...+U_sF_{[s]}=1$.
Soit $h\in \mathcal{M}(V)$ la fonction méromorphe qui 
vaut $h_j$ sur $V_{P_j}$ et $H_j=\rho(h_j)$ comme définie ci-dessus. 
\vskip 2mm
\noindent
On d\'efinit $\rho(h)$ comme le polyn\^ome 
$$
\rho(h):= H:= H_1 U_1 F_{[1]} + \cdots + H_s U_s F_{[s]}\,. 
$$
On remarque que, pour $k=0,...,d-1$,  
\begin{eqnarray*} 
{\rm Res}\, 
\left[
\begin{matrix} y^k\, H \, 
\partial_yF(y,a,b)\, dy\cr  F(y,a,b) 
\end{matrix} \right] &=&
\sum\limits_{j=1}^s {\rm Res}\, 
\left[
\begin{matrix} y^k\, H \, 
\partial_yF_j(y,a,b)\, dy\cr  F_j(y,a,b) 
\end{matrix} \right] \\
&=& 
\sum\limits_{j=1}^s {\rm Res}\, 
\left[
\begin{matrix} y^k\, H_j U_j F_{[j]} \, 
\partial_y F_j(y,a,b)\, dy\cr  F_j(y,a,b) 
\end{matrix} \right]
\end{eqnarray*} 
par les r\`egles usuelles du calcul r\'esiduel (th\'eor\`eme de dualit\'e)~; on obtient 
donc, en poursuivant suivant ces m\^emes r\`egles, et en exploitant les relations  
$U_jF_{[j]}=1-\sum_{l\ne j}U_lF_{[l]}$, $j=1,...,k$, les relations
\begin{eqnarray*} 
{\rm Res}\, \left[\begin{matrix} y^k\, H \, 
\partial_yF(y,a,b)\, dy\cr  F(y,a,b) 
\end{matrix} \right]
= \sum\limits_{j=1}^k {\rm Res}\, 
\left[
\begin{matrix} y^k\, H_j  \, 
\partial_y F_j(y,a,b)\, dy\cr  F_j(y,a,b) 
\end{matrix} \right]
\end{eqnarray*} 

soit encore

\begin{eqnarray*} 
{\rm Res}\, \left[\begin{matrix} y^k\, H \, 
\partial_yF(y,a,b)\, dy\cr  F(y,a,b) 
\end{matrix} \right]
= Tr_{V_{P_1}}y^kh_1+ \cdots +Tr_{V_{P_s}}y^kh_s=Tr_{V}y^kh
\end{eqnarray*}

Le polyn\^ome $H:=\rho(h)$ ainsi construit coïncide donc à nouveau avec l'unique polyn\^ome de degr\'e $d-1$ 
\`a coefficients dans $\mathcal{M}$ dont les coefficients  $\tau=(\tau_{0},...,\tau_{d-1})$ satisfont le syst\`eme $(S_h)$ (cette fois défini à partir des traces sur un cycle réduit et non sur un seul germe mais il reste de Cramer d'apr\`es le lemme $2$). La résolution de ce système permet donc encore de trouver les coefficients de $H$ dans le cas d'un cycle réduit. On voit ainsi (en se reportant au cas o\`u $V$ est irr\'eductible) que l'on a les trois propriétés suivantes~: 
\vskip 2mm
\noindent
\begin{itemize}
\item  
le polyn\^ome $H=\rho(h)$ v\'erifie $\partial_{a_i} H-Y\partial_{b_i} H\in (F)$ pour $i=1,...,n$~; 
\item on a, pour $j=1,...,k$, 
$$
(\rho(h))(y,a,x-ay)_{|V_j}=H_j(y,a,x-ay)_{|F_j(y,a,x-ay)=0}=h_j(x,y)~;
$$
\item $\rho(h)=0$ si et seulement si $h_{|V}\equiv 0$ ; 
en effet, $\rho(h)=0$ implique que $F_j$ divise $H_jU_jF_{[j]}$~; or $F_j$ est premier avec 
$U_j F_{[j]}$, donc $F_j$ divise $H_j$ et $H_j\equiv 0$ pour raisons de degré, soit encore $h_{j|V_{P_j}}\equiv 0$ et ce pour tout $j=1,..,s$, donc $h_{|V}\equiv 0$.
\end{itemize}
\vskip 2mm
\noindent 
Ces trois considérations montrent que 
l'application $\rho$ reste définie et bijective dans le cas de germes en plusieurs points. 
\hfill $\square$ 
\vskip 2mm
\noindent
\begin{rem}
{\rm En fait, on aurait pu définir $H$ comme étant le polynôme d'interpolation de Lagrange 
prenant les valeurs de la fonction $h(ay+b,y)$ en les racines de l'équation $f(ay+b,y)=0$~; 
on vérifie alors grâce à l'équation d'onde de choc vérifiée par les racines que $H$ ainsi construite 
vérifie les conditions $(**)$~; la construction de $\rho$ 
a l'avantage d'être plus algébrique (on ne soucie plus des singularités) et plus descriptive.}
\end{rem}
\vskip 2mm
\noindent
Le fait que l'ensemble $\{F(y,a,x-ay)=0\}$ et la fonction $H(y,a,x-ay)_{|V}$ ne dépendent pas de $a$ 
permet de choisir la direction verticale $a=0$ dans les calculs de traces. On obtient ainsi le 
moyen de reconstituer un ensemble et une forme à partir des restrictions de la trace 
d'un nombre fini de fonctions au sous-espace $a=0$ de la grassmannienne~; 
ceci paraît logique puisqu'il suffit d'une direction pour balayer un ensemble. 
Plus précisément, on a le théorème suivant~: 
\vskip 2mm
\noindent
\begin{thr}
Pour connaître un cycle $V\in\mathcal{V}$, il est nécessaire et suffisant de connaître les germes en $0$ des 
fonctions 
$$
b \to {\rm Tr}_V\, (y^k) (0,b)\,,\qquad k=0,.., {\rm Tr}_V\, (1)\,; 
$$
si $V$ est un cycle r\'eduit, pour connaître $h\in \mathcal{M}(|V|)$, il est nécessaire et suffisant de connaître les 
germes en $0$ des fonctions 
$$
b\to {\rm Tr}_V\, (y^k\, h) (0,b)\,,\quad k=0,..,{\rm Tr}_V\, (1) -1\,.   
$$
\end{thr}
\vskip 2mm
{\bf Preuve.} La preuve est contenue dans ce qui vient d'être fait. 
La connaissance de $(V,h)$ donne la connaissance des traces. A l'inverse, la connaissance des traces 
de $y^k$ sur $a=0$ permet de fabriquer le polynôme  $F(Y,0,b)$ à partir des formules de Newton 
(on peut aussi se passer de ces formules et en utilisant le syst\`eme de Cramer 
$(S)$ mais il faut alors conna{\^\i}tre les fonctions $b \to {\rm Tr}_V\, (y^k)(0,b)$ jusqu'à $k=2d-1$). Si $a=0$, le système $(S_h)$ reste non dégénéré sur le corps des germes méromorphes en $b$ d'aprés la remarque $5$ (sinon toutes les droites ``verticales'' couperaient ``mal'' $V$ ce qui est absurde). D'autre part, les traces $v_k(0,b)=Tr_Vy^kh(0,b)$ ont un sens sinon toutes les droites verticales couperaient le lieu polaire de $h$ ce qui est à nouveau absurde. Ceci permet d'obtenir le polynôme $H(Y,0,b)$ via le système $(S_h)$ en posant $a=0$ (et on montre ainsi que $H(Y,0,b)$ est bien défini). On a alors $V=\{F(y,0,x)=0\}$ et la fonction méromorphe $(x,y) \to H(y,0,x)$ co{\"\i}ncide avec $h$ sur $V$.   
\hfill $\square$ 
\vskip 2mm
\noindent
Le cas des formes de degr\'e maximal est analogue et le théorème suivant montre 
comment caractériser la trace d'une $n$-forme en termes de résidus de fonctions rationnelles 
d'une variable à coefficients méromorphes. Il suffit de l'énoncer dans le cas réduit d'après la remarque $4$. 
\vskip 2mm
\begin{thr}
Une $n$-forme $\phi$ méromorphe au voisinage de l'origine dans $\cp^{2n}_{a,b}$ 
est la trace d'une $n$-forme m\'eromorphe $\Phi$ sur un cycle r\'eduit $V$ de $\mathcal{V}$
de degr\'e vertical $d$ si et seulement s'il existe deux polynômes $F$ et $H$  tels que~: 
\begin{itemize}
\item $F\in \mathcal{U}\, [Y]$, $\deg_Y F=d$~; 
\item $H\in \mathcal{M}\, [Y]$, $\deg_Y H\leq d-1$, et $H$ v\'erifie les relations 
alg\'ebriques $\partial_{a_i}H-Y\partial_{b_i}H\in F$ pour tout $i=1,..,n$ (c'est-\`a-dire 
$H\in \mathcal{M}_F\, [Y]$)~;
\item on a l'égalité   $\phi(a,b)=Tr_V\Phi(a,b)$   avec la repr\'esentation~:  
$$
\phi(a,b)= 
{\rm Res}\, \left[
\begin{matrix} H\, (\partial_y F-\sum_j a_j\partial_{b_j} F) (y,a,b)\, dy \wedge \bigwedge\limits_{i=1}^n (db_i+y da_i) \cr
\cr 
F(y,a,b) \end{matrix}\right]\eqno(*^3) 
$$
\end{itemize} 
\end{thr}
\vskip 2mm
{\bf Preuve.} Montrons tout d'abord que si $V\in \mathcal{V}_{\rm red}$ a pour degr\'e vertical $d$ et s'\'ecrit 
$V=V_{P_1}+ \cdots + V_{P_s}$,   
les $V_{P_j}=\{f_j=0\}$, $j=1,...,s$, \'etant des germes d'hypersurfaces r\'eduites en $k$ points distincts $P_j$, 
et si $\Phi\in \mathcal{M}^n (|V|)$, alors ${\rm Tr}_V\, \Phi$ se repr\'esente sous la forme $(*^3)$. 
La donnée d'une forme méromorphe sur $|V|$ équivaut à la donnée 
d'une fonction méromorphe $h\in M(|V|)$ qui vaut $h_j$ sur $V_{P_j}$ via l'égalité $\Phi(x,y)=h(x,y)dx$. 
On pose $F_j=\Pi(V_{P_j})$, $F=\Pi(V)$, $H_j=\rho(h_j)$ et $H=\rho(h)$. D'après l'expression $(\dag\dag)$ obtenue pour la 
trace d'une $n$-forme au paragraphe $2$, les $n+1$ coefficients $w_0,...,w_{n}$ de la forme ${\rm Tr}_V\, (\Phi)$ sont 
donnés, pour $k=0,...,n$, par~:
$$
w_k(a,b):= \sum\limits_{j=1}^s 
{\rm Res}\, \left[\begin{matrix} 
\theta_j\,  y^k\, h_j\, \partial_yf(x,y)\, dx\land dy \cr   
f_j(x,y),x_1-a_1 y-b_1,\dots,x_n-a_ny-b_n 
\end{matrix}
\right]\,. 
$$
Puisque $f_j(x,y)=u_j(a,x,y)F_j(y,a,x-ay)$, où $u_j$ est inversible dans l'anneau local $\cp\{x,y-y_{P_j},a\}$, 
on peut remplacer $f_j(x,y)$ par $F_j(y,a,x-ay)$. Or on a, pour $j=1,...,s$, 
$$
\partial_y(F_j(y,a,x-ay))=\Big(\partial_y F_i-\sum\limits_{i=1}^n  a_i\partial_{b_i} F_j\Big)(y,a,x-ay)\,, 
$$
et, suivant le même argument que dans la preuve du lemme $2$, on obtient les égalités
$$
w_k(a,b)=\sum\limits_{j=1}^s 
{\rm Res}\, \left[\begin{matrix} 
y^k\, h_j(ay+b,y)\,  \Big(\partial_y F_j-\sum\limits_{i=1}^n  a_i\partial_{b_i} F_j\Big)(y,a,b)\, dy \cr 
F_j(y,a,b)\end{matrix}\right]
$$
pour $k=0,...,n$. D'apr\`es la formule $(\dag^3)$ (on l'\'ecrit pour $h=h_j$ et pour chaque facteur irr\'eductible de 
$F_j$ dans $\mathcal{M}\, [Y]$ pour se ramener au cas irr\'eductible, $F_j$ \'etant simplement suppos\'e r\'eduit ici),  
on a $(y,a,b)\to h_j(ay+b,y)-H_j(y,a,b)$ est dans l'id\'eal engendr\'e par $(y,a,b)\to F_j(y,a,b)$ dans 
$\mathcal{M}\, \{y-y_{P_j}\}$ et on a donc, pour tout $k=0,...,n$,  
$$
w_k(a,b)=\sum\limits_{j=1}^s 
{\rm Res}\, \left[\begin{matrix} 
y^k\, H_j(y,a,b)\,  \Big(\partial_y F_j-\sum\limits_{i=1}^n  a_i\partial_{b_i} F_j\Big)(y,a,b)\, dy \cr 
F_j(y,a,b)\end{matrix}\right]\,. 
$$
Il reste à raisonner comme dans la construction de 
$\rho$ (cas de plusieurs germes) pour conclure que l'on a les égalités
$$
w_k(a,b)={\rm Res}\, \left[\begin{matrix} 
y^k\, H(y,a,b)\,  \Big(\partial_y F-\sum\limits_{i=1}^n  a_i\partial_{b_i} F\Big)(y,a,b)\, dy \cr 
F(y,a,b)\end{matrix}\right]\,, \quad k=0,...,n\, 
$$
ce qui montre bien que ${\rm Tr}_V\, [h\,dx]$ est de la forme voulue. 
\vskip 2mm
\noindent
R\'eciproquement, si $\phi$ s'\'ecrit sous la forme $(*^3)$, on peut poser $V=\Pi^{-1}(F)$ et $h=\rho^{-1}(H)$ pour constater 
(en prenant les calculs \`a l'envers) que l'on a bien $\phi= {\rm Tr}_V\, (h\, dx)$. Le th\'eor\`eme $2$ est ainsi d\'emontr\'e. 
\hfill $\square$
\vskip 2mm
\noindent 
On aboutit ainsi à une nouvelle caractérisation des formes traces (autre que celle donnée 
par exemple dans \cite{H:gnus}) en termes cette fois de sommes compl\`etes de r\'esidus de 
fractions rationnelles en une variable de même dénominateur~; de telles sommes compl\`etes de 
r\'esidus s'obtiennent {\it via} l'algorithme de division euclidienne dans l'anneau $\mathcal{M}\, [Y]$.
On remarque que les coefficients de la trace s'expriment sous la forme de l'action d'op\'erateurs 
diff\'erentiels non lin\'eaires \`a coefficients polynomiaux en les coefficients $\sigma_l$ de $F$ et 
lin\'eaires en les coefficients $\tau_l$ de $H$ et les dérivées partielles des $\sigma_j$.

\section{Applications~: le th\'eor\`eme d'Abel inverse et une 
majoration de la dimension de l'espace des $q$-formes abéliennes}

\subsection{ Une forme plus forte du théorème d'Abel inverse}

On se propose dans cette section d'exploiter les r\'esultats \'etablis dans la section $3$ pour d\'emontrer
la version suivante du th\'eor\`eme d'Abel inverse~: 

\begin{thr}
Soit $V\in \mathcal{V}_{\rm red}$ un germe de cycle effectif r\'eduit intersectant proprement 
$\Delta$ (en des points $P_j$, $j=1,...,s$ distincts du point \`a l'infini) et 
de degré vertical $d$~; soit $\Phi\in {M}^n(V)$, non identiquement nulle 
sur aucune des composantes de $|V|$. Si la trace de $\Phi$ sur $V$ 
(consid\'er\'ee comme germe de $n$-forme m\'eromorphe \`a l'origine de 
$\cp^{2n}_{(a,b)}$) est le germe d'une forme rationnelle en $b$, alors  
l'ensemble analytique $|V|$ est inclus dans une hypersurface algébrique $\wt{V}$ de degr\'e $d$ telle que 
$\wt V \cap \Delta=\{P_1,...,P_s\}$ et la forme $\Phi$ se prolonge en une $n$-forme rationnelle sur 
$\pp^{n+1}(\cp)$.  
\end{thr}
\vskip 2mm
\noindent
\begin{rem} {\rm Il est essentiel de supposer que $\Phi$ n'est identiquement nulle sur aucune des composantes 
irr\'eductibles du cycle analytique $V$. 
En effet, si ce n'était pas le cas, on n'obtiendrait bien s\^ur aucune information 
relative aux composantes de $V$ sur lesquelles $\Phi$ serait identiquement nulle~; on 
pourrait seulement conclure que les autres composantes sont algébriques.}
\end{rem}
\vskip 2mm
\noindent {\bf Preuve du théorème 3.} On peut toujours supposer que $\Phi=h(x,y)dx$, où 
$h\in M(V)$. Soit $H=\rho(h)$ et $F=\Pi(V)$. On suppose que $F$ s'écrit
$$
F(Y,a,b)=Y^{d}-\sigma_{d-1}(a,b) Y^{d-1}+...+(-1)^{d-1}\sigma_0(a,b). 
$$
On introduit maintenant les germes de fonctions méromorphes \`a l'origine de 
$\cp^{2n}_{(a,b)}$~:   
$$
w_k(a,b):=
{\rm Res}\, \left[\begin{matrix} 
y^k\, H(y,a,b)\,  \Big(\partial_y F-\sum\limits_{i=1}^n  a_i\partial_{b_i} F\Big)(y,a,b)\, dy \cr 
F(y,a,b)\end{matrix}\right]\,, \, k=0,...,n\,. 
$$
D'apr\`es la formule $(\dag\dag)$ \'etablie dans la section $2$, la rationalit\'e en $b$ de 
$(a,b) \to {\rm Tr}_V\, (h\, dx)$ implique la rationalit\'e en $b$ des $w_k$, $k=0,...,n$. 
Nous aurons besoin du lemme technique suivant, concernant le comportement de la suite 
$(w_k)_{k\in \mathbb{N}}$. 
\vskip 2mm
\noindent
\begin{lem} 
La suite $(w_k)_{k\in \mathbb{N}}$ ob\'eit aux trois r\`egles suivantes~: 
\begin{itemize}
\item 
{\bf (1)} Pour tout $k\in \mathbb{N}$, pour tout $i\in \{1,...,n\}$, on a 
$$
\partial_{a_i} (w_k)=\partial_{b_i} (w_{k+1})
$$
(\'equations d'onde de choc avec décalage)~; 
\item {\bf (2)} $w_k$ est une fonction rationnelle en $b$ pour tout $k\in \mathbb{N}$~; 
\item {\bf (3)} les $w_k$ $k=0,\dots,2d-1$ vérifient le système linéaire $(\widetilde S_h)$ suivant~:
\begin{eqnarray*} 
\begin{matrix}
w_{d-1}\sigma_{d-1} &+& \cdots &+&(-1)^{d-1}w_0\sigma_0 = w_d \\
&\vdots &\vdots &\vdots & \vdots \\
w_{2d-2}\sigma_{d-1} &+&\cdots &+&(-1)^{d-1}w_{d-1}\sigma_0 = w_{2d-1} 
\end{matrix} 
\end{eqnarray*} 
\end{itemize}
\end{lem} 
\vskip 2mm
\noindent
{\bf Preuve du lemme 4.} Pour le point ${\bf (1)}$, il suffit de faire le calcul directement à partir de 
l'écriture des $w_l$ en termes de symboles r\'esiduels, telle que nous l'avons introduite 
au paragraphe $2$ (voir la formule $(\dag\dag)$)~: 
$$
w_k(a,b)={\rm Res}\, \left[\begin{matrix} 
y^k\, h\, \partial_yf(x,y)\, dx\land dy \cr   
f(x,y),x_1-a_1 y-b_1,\dots,x_n-a_ny-b_n 
\end{matrix}
\right]\,. 
$$ 
L'écriture explicite du résidu 
sous forme de repr\'esentation int\'egrale 
impliquant le noyau de Cauchy (voir par exemple 
\cite{GH:gnus}, chapitre $6$) montre que l'on peut diff\'erentier 
les symboles r\'esiduels par rapport aux param\`etres en diff\'erentiant 
sous le symbole, donc que, pour tout $k\in \mathbb{N}$, pour tout 
$i\in \{1,...,n\}$,  
\begin{eqnarray*} 
\partial_{a_i}\, (w_k) &=& - {\rm Res}\, \left[\begin{matrix} 
y^{k} \, h\, \partial_yf(x,y)\, \partial_{a_i} (L_i)\, dx\land dy \cr   
f(x,y),L_1,\dots, L_i^2,\dots, L_n 
\end{matrix}
\right] \\
&=&  {\rm Res}\, \left[\begin{matrix} 
y^{k+1} \, h\, \partial_yf(x,y)\, dx\land dy \cr   
f(x,y), L_1,\dots, L_i^2,\dots, L_n 
\end{matrix}
\right] \\
&=& -{\rm Res}\, \left[\begin{matrix} 
y^{k+1} \, h\, \partial_yf(x,y)\, \partial_{b_i} (L_i)\, dx\land dy \cr   
f(x,y), L_1,\dots, L_i^2,\dots, L_n 
\end{matrix}
\right] \\
&=& \partial_{b_i} (w_{k+1}) 
\end{eqnarray*} 
(on rappelle que $L_i(x,y,a,b):=x_i-a_i y-b_i$ pour $i=1,...,n$). 
\vskip 2mm
\noindent
On montre le point ${\bf (2)}$ par r\'ecurrence sur $k$. C'est vrai pour $k=0,\dots,n$ par hypothèse. 
On suppose la propri\'et\'e vraie jusqu'au rang $k-1$. 
\vskip 1mm
\noindent
On suppose dans un premier temps que $w_k\ne 0$ pour tout $k\in \mathbb{N}$. 
On va reprendre l'astuce utilisée par G. Henkin et M. Passare 
\cite{hp:gnus} pour éviter les pôles simples (ce point, qui est aussi directement 
inspir\'e de l'article fondateur d'Abel \cite{ab:gnus, B:gnus}, para{\^i}t \^etre 
le point clef de la d\'emonstration).  
On pose pour cela $w_k'=w_k+cw_{k-1}$ où $c\in\cp^*$. D'après le point 
{\bf (1)}, on a, pour $i=1,...,n$, les deux égalités
\begin{eqnarray*}
\partial_{b_i} (w_k')&=& \partial_{b_i}(w_k+cw_{k-1})= \partial_{a_i}(w_{k-1}+c w_{k-2})=\partial_{a_i} (w_{k-1}') \\
\partial_{b_i} (w_k')&=& (\partial_{a_i}+c\partial_{b_i}) (w_{k-1})\,. 
\end{eqnarray*}
Ainsi, la fonction $\partial_{b_i} (w_k')$ admet-elle 
les deux primitives distinctes $w_{k-1}$ et $w_{k-1}'$ dans les deux directions linéairement indépendantes 
$a_i$ et $a_i+cb_i$ (car $w_{k-2}\ne 0$ par l'hypothèse faite pour l'instant).
Or par hypothèse $\partial_{b_i} (w_k')=(\partial_{a_i}+c\partial_{b_i})(w_{k-1})$ est rationnelle en $b_i$~;  
le fait que cette fonction méromorphe ait deux primitives distinctes dans deux directions linéairement 
indépendantes exclut que dans la décomposition en éléments simples de cette fraction rationnelle 
dans $\overline{\mathcal{M}_{a,\hat b_i}}(b_i)$ (ici $\overline{\mathcal{M}_{a,\hat b_i}}$ d\'esigne 
une cl\^oture int\'egrale du corps des germes de fonctions 
m\'eromorphes en $a,b_1,...,\widehat{b_i},...,b_n$ \`a l'origine de $\cp^{2n-1}$), on puisse rencontrer un terme 
du type 
$$
\frac{\alpha_{a,\hat b_i}}{b_i-\beta(a,\hat b_i)}\,, 
$$
ce qui permet d'intégrer selon la variable $b_i$ et d'obtenir ainsi 
une nouvelle fraction rationnelle en $b_i$ (il ne peut plus y avoir de pôles simples et 
donc de logarithme dans la primitive). Puisque $w_k=w_l'-cw_{k-1}$, alors $w_k$ est rationnelle en $b_i$ 
puisque \`a la fois $w_l'$ et  $w_{l-1}$ le sont.  
\vskip 1mm
\noindent
Si maintenant $w_k=0$ pour un $k\in \mathbb{N}$ (on prend le plus petit $k$ pour lequel ceci se produit), 
alors, puisque $0=\partial_{a_i} (w_k)=\partial_{b_i} (w_{k+1})$ pour $i=1,...,n$, 
la fonction $w_{k+1}$ ne dépend pas de $b$~; comme $\partial_{b_i} (w_{k+2})=\partial_{b_i} (w_{k+1})$, 
la fonction $w_{k+2}$ est en fait polynômiale de degré $1$ en $b$~; de proche en proche, on montre alors que 
$w_{k+j}$ est polynômiale de degré $\le j$ en $b$.
\vskip 1mm
\noindent
Dans l'une ou l'autre situation envisag\'ees, nous avons prouv\'e le point {\bf (2)}. 
\vskip 2mm
\noindent
Le point {\bf (3)} r\'esulte simplement des identit\'es 
$$
{\rm Res}\, \left[\begin{matrix} 
y^k\, H(y,a,b)\, F(y,a,b)\,  \Big(\partial_y F-\sum\limits_{i=1}^n  a_i\partial_{b_i} F\Big)(y,a,b)\, dy \cr 
F(y,a,b)\end{matrix}\right]=0 
$$
pour tout $k\in \mathbb{N}$ (propri\'et\'e du calcul r\'esiduel), qui, une fois $F$ d\'evelopp\'e au num\'erateur, font 
bien appara{\^\i}tre que les $w_k$ v\'erifient le système linéaire $(\widetilde S_h)$.  
\vskip 2mm
\noindent
Le lemme $4$ est donc ainsi démontré. \hfill $\square$ 
\vskip 2mm
\noindent
\begin{rem} 
{\rm On retrouve ici, en prouvant la clause {\bf (1)} du lemme $4$, 
que la trace d'une $n$-forme est une forme fermée (en dehors de son lieu polaire) 
en utilisant son expression $(\dag\dag)$ en termes des $w_k$ introduits ci-dessus~; ceci 
est en fait une conséquence du fait que
$$
d(\Phi\land[V])=d\Phi\land [V]=0
$$
dans le cas où $\Phi$ est une forme de degré maximal sur $V$ et du fait que $d$ commute avec le 
pull-back et l'image directe.}
\end{rem}
\vskip 2mm
\noindent
Avant de reprendre la preuve du théoréme $3$, prouvons le second lemme auxiliaire suivant~: 
\vskip 2mm
\noindent
\begin{lem}
Le système $(\widetilde S_h)$ est un système de Cramer. 
\end{lem}
\vskip 2mm
\noindent
\noindent {\bf Preuve du lemme 5.} On raisonne par l'absurde en supposant que $(\widetilde S_h)$ 
est un syst\`eme dégénéré. On peut donc, dans ce cas, 
construire un polynôme $Q\in \mathcal{M}[Y]$, unitaire et de degré $d$, différent de $F$, et tel que 
$$
{\rm Res}\, \left[\begin{matrix} 
y^l\, H(y,a,b) Q(y,a,b)\, \Big(\partial_y F-\sum\limits_{i=1}^n  a_i\partial_{b_i} F\Big)(y,a,b)\, dy \cr 
F(y,a,b)\end{matrix}\right]=0 
$$
pour $k=0,...,d-1$. Le théorème de dualité dans le contexte 
algébrique (voir le rappel en fin de section $2$) implique que le polynôme 
$Q\, H\, (\partial_yF-\sum_ia_i\partial_{b_i}F)$ est dans l'idéal engendré 
par $F$ dans $\mathcal{M}[Y]$. Puique $Q\ne F$ et que $F$ est réduit, il existe forcément un facteur irréductible  
$F_j$ de $F$ (une fois $F$ décomposé dans $\mathcal{M}[Y]$) qui ne divise pas $Q$. 
Dans ce cas, le polynôme $F_j$ divise $H(\partial_yF-\sum_ia_i\partial_{b_i}F)$. Deux situations 
sont alors envisageables~: 
\vskip 2mm
\noindent
\begin{itemize}
\item le polynôme $F_j$ divise $H$, mais alors $H(y,0,x)_{|V_j=\{F_j(y,0,x)=0\}}=0$ 
et on tombe sur le cas pathologique (et exclus) où la forme $\Phi$ est identiquement nulle sur 
la composante $V_j$~; 
\item le polynôme $F_j$ divise $\partial_yF-\sum_ia_i\partial_{b_i}F$~; en notant $F=F_j\wt{F}$ 
($F_j$ et $\wt{F}$ étant premiers entre eux), on a alors~:  
$$
\partial_yF-\sum\limits_{i=1}^n a_i\partial_{b_i} F =F_j\wt{F}\, \Big(\partial_yF_j-
\sum\limits_{i=1}^n a_i\partial_{b_i}F_j\Big)+F_j\, \Big(\partial_y\wt{F}-\sum\limits_{i=1}^n a_i\partial_{b_i}\wt{F}\Big) \,, 
$$
ce qui implique que $F_j$ divise $\wt{F}(\partial_yF_j-\sum_ia_i\partial_{b_i}F_j)$ et, par 
le lemme de Gauss, que $F_j$ divise $\partial_yF_j-\sum_ia_i\partial_{b_i}F_j$~; pour des raisons de degré, 
la seule possibilité est $\partial_yF_j-\sum_ia_i\partial_{b_i}F_j=0$ ce qui implique que $\partial_y(F_j(y,a,x-ay))=0$~; 
mais alors l'équation de la branche $V_j=\Pi(F_j)$ ne dépend donc pas de $y$ et $V_j$ est forcément la droite verticale $x=0$, 
ce qui est exclus puisque $V$ intersecte proprement cette droite. 
\end{itemize} 
\vskip 2mm
\noindent
Les deux situations pouvant poser problème ayant de fait été exclues par hypothèses, le 
système $(\widetilde S_h)$ est bien de Cramer, ce qui achève la preuve du lemme $5$. \hfill $\square$ 
\vskip 2mm
\noindent
{\bf Suite de la preuve du théorème $3$.} En combinant les lemmes $4$ et $5$, on voit que 
$\sigma_0,...,\sigma_{d-1}$ s'expriment rationnellement en fonction des $w_k$, $k=0,...,2d-1$.  
Les coefficients de $F$ sont donc rationnels en $b$ et 
l'ensemble 
$$
\{(x,y)\in\cp^2 \,;\, F(y,0,x)=0\} 
$$
définit une hypersurface algébrique de $\mathbb{P}^{n+1}$ qui contient $V$ et dont 
le degré est donné par 
$$
Tr_{\wt{V}} (1)= {\rm Res}\, \left[\begin{matrix} \partial_y F(y,a,b) \, dy \cr F(y,a,b)\end{matrix}\right] =d\,, 
$$
ce qui montre le premier volet de la conclusion du théorème (il existe une hypersurface algébrique de degré 
$d$ contenant le support de $V$ et ne rencontrant la droite projective $\Delta$ qu'aux points $P_1,...,P_s$).
\vskip 2mm
\noindent  
On peut d'ailleurs faire à ce stade de la preuve la remarque suivante~:  
\vskip 2mm
\noindent
\begin{rem} {\rm Si $F\in \mathcal{U}\, [Y]$, avec de plus des coefficients rationnels en $b$ (ce qui est le cas ici) 
on sait, d'après la remarque 6, que le lieu polaire des coefficients de $F$, considérés comme des 
éléments de $\mathcal{M}_a (b)$, o\`u $\mathcal{M}_a$ désigne le corps des germes de fonctions méromorphes 
en $a$ à l'origine de $\cp^n$, ne saurait dépendre de $b$~; les coefficients d'un tel polynôme $F$ 
sont donc automatiquement polynômiaux en $b$ dès qu'ils sont rationnels en $b$~; en fait, on montre aisément 
en écrivant les conditions $(*)$ qu'une condition nécessaire et suffisante pour qu'un élément $Y^d-\sigma_{d-1} Y +\cdots+(-1)^{d-1}\sigma_0$ de 
$\mathcal{U}\, [Y]$ ait ses coefficients rationnels en $b$ est que $\sigma_{d-1}$ soit une fonction affine de $b_i$, 
$i=1,...,n$. Dans notre contexte, ce coefficient correspond à la trace de $y$ sur $V$ et 
on retrouve ainsi le théorème de J. Wood \cite{W:gnus}~: 
une condition nécessaire et suffisante pour l'existence de $\wt{V}$ est le fait que ${\rm Tr}_V (y)$ 
soit affine en $b$. Construire la trame logique reliant les énoncés d'Abel inverse et de Wood était en fait la première 
motivation de ce travail, ce au regard de l'effectivité de la construction de J. Wood.} 
\end{rem}
\vskip 2mm
\noindent
{\bf Suite et fin de la preuve du théorème $3$.} On montre maintenant la 
rationalité de $\Phi$. On introduit pour cela les coefficients holomorphes 
$\xi_k$, $k\in \mathbb{N}$, impliqués dans la définition des traces ${\rm Tr}_V\, (y^k\, dx)$, 
$k\in \mathbb{N}$ et définis par 
$$
\xi_k:= {\rm Res}\, \left[\begin{matrix} 
y^k\, \Big(\partial_y F-\sum\limits_{i=1}^n  a_i\partial_{b_i} F\Big)(y,a,b)\, dy \cr 
F(y,a,b)\end{matrix}\right]\,,\quad k\in \mathbb{N}\,. 
$$
Les coefficients de $H$ v\'erifient (on le voit en développant $H$) le système
\begin{eqnarray*} 
(\check S_h) \qquad \begin{matrix} \xi_{d-1}\tau_{d-1} &+& \cdots &+&\xi_0\tau_0 = w_0 \\
&\vdots & \vdots &\vdots & \vdots \\
\xi_{2d-2}\tau_{d-1} &+& \cdots &+&\xi_{d-1}\tau_0 = w_{d-1} 
\end{matrix} 
\end{eqnarray*} 
Comme précédemment (voir la preuve du lemme $5$), 
on peut montrer que ce système est un système de Cramer. 
Puisque les coefficients de $F$ sont rationnels en $b$, les fonctions $\xi_l$ le sont et les coefficients $\tau_k$ 
le sont aussi ce qui permet de conclure que $H$ a ses coefficients rationnels en $b$. 
D'aprés la preuve du théorème $1$, la forme $H(y,0,x)\, dx$ est bien définie. Elle définit donc une forme $\wt{\Phi}$ rationnelle sur $\pp^{n+1}(\cp)$ prolongeant la forme 
$\Phi$. Ceci achève la preuve du théorème $3$. \hfill $\square$ 
\vskip 2mm
\noindent
\begin{rem}
{\rm 

a) Cette preuve s'adapte à la version locale d'Abel-inverse: si la trace se prolonge à un ouvert $\wt{D}^*$ contenant $D^*$, alors l'ensemble $V$ et la forme $\Phi$ se prolongent à l'ouvert $1$-concave $\wt{D}$ (contenant $D$) dont le dual est $\wt{D}^*$.

b) Les bijections $\Pi$ et $\rho$ permettent de montrer que toute propriété vérifiée par 
${\rm Tr}_V\, (h(x,y)dx)$ se répercute en une propriété (en les variables 
$b_1,...,b_n$) pour les coefficients des polynômes $F=\Pi(V)$ et $H=\rho(V)$. 
C'est alors le fait que les applications $V=\Pi^{-1}(F)$ et $h=\rho^{-1}(H)$ ne se soucient pas du 
comportement en $a$ de $F$ et $H$ qui permet de basculer les propriétés de ${\rm Tr}_V \, (h(x,y)dx)$ en des  
propriétés relatives à $V$ et $h$~; par exemple, si la trace de $\Phi$ sur $V$ est algébrique, on doit pouvoir 
montrer (toujours avec l'argument des deux primitives) que $F$ et $H$ sont algébriques en $b$~; 
Ainsi en multipliant alors les ``conjugués'' de $F$ (définis {\it via} les automorphismes de l'extension galoisienne 
finie de $\cp\{a\}(b)$ engendrée par les coefficients $\sigma_k(a,b)$ de $F$), on va fabriquer deux nouveaux 
polynômes $\wt{F}$ et $\wt{H}$ à coefficients rationnels vérifiant les conditions $(*)$ et $(**)$~; ils donneront 
alors naissance {\it via} $\Pi^{-1}$ et $\rho^{-1}$ à un ensemble algébrique $\wt{V}$ (dont le degré sera le degré de $F$ 
multiplié par le degré de l'extension) et à une forme rationnelle $\wt{\Phi}$ tels que $V\subset\wt{V} $ et 
$\wt{\Phi}_{|V}=\Phi$. On retrouve dans ce cas le théorème d'Abel inverse version algébrique montré par 
S. Collion (voir \cite{C:gnus}).}
\end{rem}
\vskip 2mm
\noindent
Il semble également intéressant de souligner qu'une fois le degré vertical $d$ précisé, le théorème 
d'Abel inverse permet d'affirmer que toutes les solutions $(\sigma_0,...,\sigma_{d-1},\tau_0,...,\tau_{d-1})
\in (\mathcal{O}_{(a,b)})^{d}\times (\mathcal{M}_{(a,b)})^d$ 
du système linéaire du premier ordre en les inconnues $\tau_j$, $j=0,...,d-1$, 
différentiel polynomial du premier ordre en les inconnues $\sigma_j$, $j=0,...,d-1$, avec second membre donné $\Psi$,  
s'écrivant~: 
\begin{eqnarray*} 
F &=& Y^d-\sigma_{d-1} Y^{d-1}+ \cdots +(-1)^d \sigma_0 \in \mathcal{U}[Y]   \\   
H &=& \tau_{d-1} Y^{d-1} + \cdots + \tau_1 Y + \tau_0 \in \mathcal{M}_F\, [Y] \\ 
&{\rm Res}& \, \left[
\begin{matrix} H\, (\partial_y F-\sum_j a_j\partial_{b_j} F) (y,a,b)\, dy \wedge \bigwedge\limits_{i=1}^n (db_i+y da_i) \cr
\cr 
F(y,a,b) \end{matrix}\right] = \Psi(a,b) 
\end{eqnarray*} 
sont des solutions rationnelles (si $F$ et $H$ sont supposés premiers entre eux) 
ou tout au moins algébriques (si cette dernière restriction n'est pas imposée) 
dès que le second membre (c'est-à-dire la $n$-forme $\Psi$) est une forme 
rationnelle. Comme le lieu polaire des $\sigma_l\in {\cal M}_{(a,b)}$ ne d\'epend pas de $b$ d\`es 
que $F=Y^d-\sigma_{d-1} Y^{d-1}+ \cdots +(-1)^d \sigma_0$ v\'erifie les conditions $(*)$ (voir la remarque $6$), on voit 
d'ailleurs que la rationalit\'e du second membre $\Psi$ implique la rationalit\'e 
des solutions $(\sigma,\tau)$ du syst\`eme dans $(\mathcal{M}_{(a,b)})^{2d}$ (on perturbe 
$a=0$ en prenant $a$ g\'en\'erique voisin de $0$ de mani\`ere \`a se ramener au cas o\`u 
$F\in \mathcal{U}\, [Y]$). Cette remarque indique 
que l'on peut concevoir le théorème Abel inverse comme un résultat de {\it rigidité} 
relatif à un système différentiel non linéaire d'un type très particulier 
(linéaire d'ordre $0$ en $\tau$, d'ordre $1$ en les d\'eriv\'ees de 
$\sigma$, polynômial en $\sigma$). Il nous para{\^\i}t dès lors important de 
formuler de manière identique (c'est-\`a-dire en terme de rigidité d'un certain système 
différentiel du même type) les théorèmes du type Abel inverse o\`u la grassmannienne 
se trouve remplacée par une famille de courbes (ou bien la variété ambiante 
$\pp^{n+1}(\cp)$ rempla\c cée par exemple par une variété torique complète 
simpliciale de dimension $n+1$ correspondant à un polytope de Delzant et par conséquent 
plongée dans un certain $\pp^N(\cp)$)~; on trouve par exemple une généralisation du théorème 
classique d'Abel inverse dans \cite{F_1:gnus}, pouvant servir de point de départ pour ces nouvelles 
approches.   
\vskip 2mm
\noindent 
Notons aussi que, si seulement la machinerie élémentaire du calcul résiduel 
en une variable s'est trouvée impliquée ici, il est vraisemblable que les généralisations 
évoquées nécessiteront une utilisation des outils du calcul de résidus en plusieurs 
variables complexes cette fois dans toute leur puissance. 
\vskip 2mm
\noindent
J'espère revenir sur ces diverses questions et généralisations dans des travaux ultérieurs.   

\subsection{A propos de la dimension de l'espace des $q$-formes abéliennes}

Soit $V$ une hypersurface algébrique réduite de $\pp^{n+1}$.
On cherche à caractériser l'ensemble $\om^n(V)$ des formes abéliennes sur $V$, {\it i.e.} des formes dont la trace sur $V$ est nulle. (Car holomorphe sur toute la grassmanienne). On peut toujours choisir un système de coordonnées pour lequel $V$ est coupée proprement par $\Delta=\{x=0\}$ et ce en $s$ points de l'espace affine $\cp^{n+1}$. Au voisinage de cette droite, $V$ permet donc de définir un élément $V_1+\cdots +V_s$ de $\mathcal{V}_{\rm red}$ de degré vertical $d$. D'aprés le théorème 1, toute forme $\Phi=hdx$ méromorphe 
sur $V$ est uniquement déterminée par les $d$ fonctions de $v_k(0,b)=Tr_V y^kh(0,b) ,k=0,..,d-1$. Pour $a=0$ on remarque que l'on a en fait $v_k(0,b)=w_k(0,b)$. Or les $w_k\,,k=0,..,n$ sont les coefficients de la trace, et dans le cas des formes de trace nulle, on a $w_0=...=w_n=0$ et $w_{n+j}$ est polynômial en $b$ de degré $\le j-1$ pour tout $j$ (d'après la preuve de la propriété $3$ du lemme $4$).
Une forme abélienne est ainsi uniquement déterminée par la collection de $d-n-1$ polynômes de $n$ variables $P_i(b):=w_{n+i}(0,b),\,i=1,..,d-n-1$ avec $\deg P_i\leq i-1$ 
(il n'y a donc pas de forme abélienne non nulle dans le cas $d<n$). Ainsi l'espace des $n$ formes 
abéliennes sur une hypersurface algébrique réduite $V\subset\pp^{n+1}(\cp)$ de degré $d$ s'identifie à un sous-espace de l'espace vectoriel des 
$d-n-1$ uplets de polynômes $(P_1,...,P_{d-n-1})$ tels que $\deg P_i\leq i-1$~; comme ce dernier 
espace est de dimension $\Big(\begin{matrix} d-1 \cr n+1 \end{matrix}\Big)$, la dimension de l'espace 
vectoriel $\omega^n(V)$ est bien majorée par $\Big(\begin{matrix} d-1 \cr n+1 \end{matrix}\Big)$. 
\vskip 2mm
\noindent
Pour montrer que la majoration est atteinte, on considère les formes rationnelles s'écrivant localement au voisinage de la droite verticale sous la forme:
$$
\Phi(x,y)=\frac{P(x,y)}{\partial_yf(x,y)}dx
$$
où
$$
P(x,y)=a_0y^{d-n-2}+a_1(x)y^{d-n-3}+...+a_{d-n-2}(x)\quad deg\;a_i\le i
$$
est un polynôme en $(x,y)$ de degré total inférieur ou égal à $d-n-2$ et $f$ est un polynôme de degré $d$ (non divisible par $x$, comme dans le cas des germes) donnant l'équation affine $f(x,y)=0$ de $V$ au voisinage de la droite verticale $x=0$. Dans ce cas, on a:
$$
Tr_V\Phi(a,b)=\sum_{k=0}^n  \Big(w_k(a,b)\big(\sum_{|I|=k}da_I\land db_{I^c}\big)\Big)
$$
et à partir de l'écriture $(\dag\dag)$ du paragraphe 2, puis en utilisant le lemme 2, on obtient cette fois
$$
w_k(a,b)={\rm Res}\, \left[\begin{matrix} 
y^k\, P(ay+b,y)\, dy \cr 
F(y,a,b)\end{matrix}\right]. 
$$
Or, si $k\le n$, $deg(y^k\, P(ay+b,y))\le d=deg\,F-2$  d'où la nullité des $w_k(a,b)$ pour tout $k=0,...,n$. Toute forme ainsi définie est donc une forme abélienne. De plus, toujours par le théorème de dualité, on constate que le polynôme en $y$ $P(ay+b,y)$ (et donc le polynôme $P\in\cp[x,y]$) est entièrement caractérisé par les $w_k(a,b)\;k=n+1,....,d-1$ . Ces dernières fonctions caractérisant entièrement la forme de trace nulle $\Phi_{|V}$, la dimension de l'espace vectoriel des formes abéliennes se trouve minorée par la dimension de l'espace vectoriel des polynômes $P(x,y)$ de degré inférieur ou égal à $d-n-2$ qui est exactement $\Big(\begin{matrix} d-1 \cr n+1 \end{matrix}\Big)$. On retrouve ainsi l'égalité:
$$
dim\;\omega^n(V)=\Big(\begin{matrix} d-1 \cr n+1 \end{matrix}\Big)
$$
avec une écriture affine explicite des formes abéliennes de degré maximales sur une hypersurface algébrique. On peut noter que travailler par dualité avec la trace permet de ne pas se soucier du comportement de $\Phi$ à l'infini (contrairement aux caractérisations des formes abéliennes en général, par exemple \cite{he:gnus}). Le fait que $\Phi$ (dans son écriture affine) n'ait pas de pôles sur l'hyperplan à l'infini (sauf éventuellemnt sur $V$) est impliqué par la nullité de la trace: on évite ainsi les changements de carte de $\pp^{n+1}$. 
\vskip 2mm
\noindent
On s'intéresse maintenant à l'espace $\om^q(V)$ des $q$-formes abéliennes sur une hypersurface algébrique 
réduite $V\subset\mathbb{P}^{n+1}$ de degré vertical $d$. 
Au voisinage de la droite verticale, toute $q$-forme rationnelle $\Phi$ peut s'écrire
$$
\Phi=\sum_{|I|=q}h_Idx_I\,, 
$$
où les $h_I$ sont des fonctions rationnelles de $(x,y)$. D'après le théorème 1, pour connaître $h_I$ 
il est nécessaire est suffisant de connaître les $d$ fonctions 
$$
t_{I,k}(0,b):={\rm Tr}_V\, [y^k h_I] (0,b)\,,\qquad k=0,\dots,d-1\,, 
$$
définies par
$$
t_{I,k}={\rm Res}\, \left[\begin{matrix}Y^kh_I(Y,b)\partial_YF(Y,0,b)dY \cr F \end{matrix}\right]\,,   
$$
où $F=\Pi(V)$ (on se restreint ici à la carte affine $\cp^{n+1}$ dans laquelle on travaille et 
dans laquelle $V=\{f(x,y)=0\}$). Si l'on se restreint aux droites verticales, l'expression des coefficients de la trace d'une 
$q$-forme se simplifie~; plus précisément, on a défini dans l'introduction le courant $T$ par
$$
T=\Phi\land df\land \Big(\bigwedge_{j=1}^n dL_j\Big)  
\land \bar{\partial}\Big(\frac{1}{f}\Big)\land \Big(\bigwedge_{j=1}^n {\partial}\Big(\frac{1}{L_i}\Big)\Big)\,.  
$$
On a alors ici
$$
\Phi\land df \land 
\Bigg(\bigwedge_{j=1}^n dL_j\Bigg) 
=\sum_{|I|=q}h_I\partial_yfdx_I\land dy\land \bigwedge_{j=1}^n (dx_j-a_jdy-yda_j-db_j)\,.  
$$
On peut supposer ici qu'il n'y a pas les termes $a_idy$ puisqu'on s'intéresse aux coefficients de 
la trace restreints à $a_1=\cdots= a_n=0$ pour caractériser les $h_i$~; de plus on ne garde que les termes 
$da_K\land db_L$ pour $|K|+|L|=q$ puisque l'on fait agir $T$ sur les formes-test de $D\times D^*$ de 
bidegré $(2n-q,2n)$ (en $(a,b)$). Dans ce cas, on trouve une expression simplifiée
$$
\Phi\land df\land \bigwedge_{j=1}^n dL_j =
\sum_{|I|=q} h_I\, \partial_yf\, dx_I\land dx_{I^c}\land dy\land
\Big(\sum\limits_{J\subset I,|J|=0}^q \pm y^{|J|}da_J\land db_{I\setminus J} \Big)
$$
et l'on a
$$
{\rm Tr}_V\,  [\Phi](0,b)=\sum_{|I|=q}\sum_{J\subset I,|J|=0}^q 
\pm {\rm Res}\, \left[\begin{matrix}
Y^{|J|}h_I(Y,b)\partial_YF(Y,0,b)dY \cr F \end{matrix}\right]  \, da_J\land db_{I\setminus J}\,
$$
Comme dans le cas des $n$-formes, le fait que la trace soit nulle implique alors 
$t_{I,k}=0$ pour tout $k=0,...,q$ et $\deg_{b} [t_{I,q+k}] <k$ pour tout $k>0$, et ce 
pour tout multi-indice $I\subset\{1,..,n\}$ de longueur $q$. Chaque $h_I$ est alors uniquement 
déterminé par la collection de polynômes 
$$
\{t_{I,q+k}(0,b)\in\cp[b]\, ;\, k=1,\dots,d-q-1\,;\, \deg_b [t_{I,q+k}]<k\}\,.
$$
Puisqu'il y a un nombre $\Big(\begin{matrix} n\cr q \end{matrix}\Big)$ de multi-indice $I\subset\{1,..,n\}$ de longueur $q$ et que l'ensemble des collections possibles pour chaque $I$ est un espace vectoriel sur $\cp$ de dimension 
$\Big(\begin{matrix} d+n-q-1 \cr n+1 \end{matrix}\Big)$
on retrouve le fait que l'espace des $q$-formes abéliennes sur une hypersurface algébrique de degré $d$ est 
un espace vectoriel de dimension finie majorée par le nombre de Castelnuovo $\pi_q(d,2,n)$~:
$$
\dim \om^q(V)\le \pi_q(d,2,n):=\Big(\begin{matrix} n 
\cr q \end{matrix}\Big)\,  \Big(\begin{matrix} d+n-q-1 \cr n+1 \end{matrix}\Big)
$$
Cependant, contrairement \`a ce qui se passe dans le cadre particulier des formes abéliennes de degré maximal 
(pour lesquelles on peut également profiter du fait que la trace soit $d$-fermée), on ne peut pas conclure 
dans ce cas que la borne est atteinte. En fait si on suppose chaque $h_I$ sous la forme $\frac{P_I}{\partial_yf}$ avec $deg P_I\le d-q-2$ on obtient des formes abéliennes, mais cette fois deux écritures distinctes peuvent être égales modulo $(f,df)$, donc donner la même forme sur $V$. Ainsi, le défaut entre la borne de Castelnuovo et la dimension de $\omega^q(V)$ dépend cette-fois de la différentielle de $f$ donc des singularités de $V$: par exemple, si $V$ est une hypersurface lisse de $\pp^3$, l'espace vectoriel $\om^1(V)$ est réduit à $0$ ce qui n'est à priori pas le cas si $V$ est singulière. Le cas des $q$-formes abéliennes $(q<n)$ est beaucoup plus délicat que celui des formes de degré maximales (voir par exemple \cite{he:gnus}) et est peut-être un bon terrain de recherche.

\end{document}